\documentclass[amscd]{amsart}
\usepackage{amsmath,amssymb,amscd}
\usepackage{xy}
\xyoption{all}

\newcommand{\A}{\ensuremath{\mathcal A}}
\newcommand{\B}{\ensuremath{\mathcal B}}
\newcommand{\C}{\ensuremath{\mathbb C}}
\newcommand{\CC}{\ensuremath{\mathcal C}}

\newcommand{\E}{\ensuremath{\mathcal E}}

\newcommand{\K}{\ensuremath{\mathcal K}}

\newcommand{\R}{\ensuremath{\mathbb R}}
\renewcommand{\S}{\ensuremath{\mathcal S}}

\newcommand{\Z}{\ensuremath{\mathbb Z}}

\newcommand{\Cliff}{{\operatorname{Cliff}}}

\newcommand{\cs}{{C^{\star}}}

\newcommand{\into}{\hookrightarrow}

\renewcommand{\epsilon}{{\varepsilon}}
\newcommand{\ev}{\operatorname{ev}}

\newcommand{\ga}{{\Gamma}}

\newcommand{\gtimes}{{\hat{\otimes}}}
\newcommand{\ideal}{\vartriangleleft}
\newcommand{\id}{\operatorname{id}}
\newcommand{\Ind}{{\operatorname{Ind}}}

\newcommand{\od}{\operatorname{od}}

\newcommand{\rank}{{\operatorname{rank}}}
\newcommand{\range}{{\operatorname{range}}}

\newcommand{\spinc}{{spin}$^c$\!}

\newcommand{\vlim}[1]{\underset{#1}{\varinjlim}\ }

\newcommand{\Aut}{\operatorname{Aut}}
\newcommand{\FB}{\mathcal{FB}}
\newcommand{\EE}{{\mathfrak E}}
\newcommand{\FF}{{\mathfrak F}}
\renewcommand{\ga}{g}
\renewcommand{\gg}{\tilde{g}}

\newcommand{\metric}{{\langle \cdot, \cdot \rangle}}
\newcommand{\minus}{\backslash}
\newcommand{\Prim}{\operatorname{Prim}}
\newcommand{\Res}{\operatorname{Res}}

\newtheorem{thm}{Theorem}[section] 
\newtheorem{dfn}[thm]{Definition}

\newtheorem{cor}[thm]{Corollary}

\newtheorem{ex}[thm]{Example}
\newtheorem{lem}[thm]{Lemma}

\newcommand{\Pf}{{\em Proof}. }
\newcommand{\Pff}{{\em Proof} }
\newtheorem{prop}[thm]{Proposition}

\newtheorem{rmk}[thm]{Remark}
\newtheorem{thm*}{Theorem}
\newcommand{\EPf}{\hfill \qed}
  {\addvspace{\bigskipamount}\noindent{\bf Example\ }}%
  {}

\newtheorem{ThomThm}[thm]{Thom Isomorphism Theorem}

\newenvironment{tightlist}%
{\begin{list}{$\bullet$}{\parsep=0pt \itemsep=0pt 
                         \topsep=6pt \partopsep=\parskip}}%
{\end{list}}

\newcounter{ictr}
\newenvironment{ilist}{\begin{list}
                         {(\roman{ictr})}
                         {\usecounter{ictr}
                          \parsep=0pt \itemsep=0pt 
                          \topsep=6pt \partopsep=\parskip
                          \setlength{\leftmargin}{0.5truein}}}
                      {\end{list}}

\newcounter{nctr}

\begin{document}%
\title{A Thom Isomorphism for Infinite Rank Euclidean Bundles}          
\author{Jody Trout}             
\email{jody.trout@dartmouth.edu}       
\address{6188 Bradley Hall\\ 
         Dartmouth College\\
         Hanover, NH 03755}
\thanks{The author was partially supported by NSF Grants DMS-9706767 and DMS-0071120}
\subjclass{19, 46, 47, 55, 57, 58}
\begin{abstract}
An equivariant Thom isomorphism theorem in operator $K$-theory is formulated and proven for infinite rank Euclidean vector bundles over finite dimensional Riemannian manifolds. The main ingredient in the argument is the construction of a non-commutative $\cs$-algebra associated to a  bundle $\EE \to M$, equipped with a compatible connection $\nabla$, which plays the role of the algebra of functions on the infinite dimensional total space $\EE$. If the base $M$ is a point, we obtain the Bott periodicity isomorphism theorem of Higson-Kasparov-Trout \cite{HKT98} for infinite dimensional Euclidean spaces. The construction applied to an even {\it finite rank} \spinc-bundle over an even-dimensional proper \spinc-manifold reduces to the classical Thom isomorphism in topological $K$-theory. The techniques involve non-commutative geometric functional analysis.
\end{abstract}



\maketitle

\section{Introduction}

First, we review the classical Thom isomorphism in topological $K$-theory, in a smooth version appropriate for our infinite dimensional
generalization. Let $M$ be a smooth finite dimensional manifold on which the compact group $G$ acts via diffeomorphisms. The
equivariant topological $K$-theory group of $M$ (with compact supports) can be defined \cite{Seg68} as the abelian group $K^0_G(M)$
generated by $G$-homotopy equivalence classes $[\sigma]$ of smooth equivariant morphisms $\sigma : F_1 \to F_2$ of smooth complex vector
$G$-bundles (with finite rank) which are isomorphisms off a compact subset of $M$. Two such bundle morphisms
$\sigma$ and $\sigma'$ are $G$-homotopic if there is a $G$-bundle morphism $\Sigma$ on $M \times [0,1]$, also with compact support, such
that $\Sigma|_{M \times \{0\}} \cong \sigma$ and $\Sigma|_{M \times \{1\}} \cong \sigma'$. Addition in $K^0_G(M)$ is given by direct sum:
$[\sigma] + [\sigma'] = [\sigma \oplus \sigma'].$ The map $M \mapsto K^0_G(M)$ defines a contravariant functor from the category of smooth
$G$-manifolds (and smooth proper $G$-maps) to the category of abelian groups.

If $p : E \to M$ is a smooth vector $G$-bundle of finite rank on $M$ (either real or complex), then
$K^0_G(E)$ has a right $K^0_G(M)$-module structure via the composition
$$\xymatrix{  K^0_G(E) \times K^0_G(M) \ar[r]^-{\boxtimes} \ar@/_1pc/[rr]_-{\otimes} & K^0_G(E \times M) \ar[r]^-{(\id \times p)^*} & K^0_G(E)}$$ 
where the first homomorphism is given by the outer tensor product $\boxtimes$ and the second homomorphism is induced by the
proper $G$-map $\id \times p : E \to E \times M$. Give the vector bundle $E$ a $G$-invariant bundle metric $(\cdot , \cdot)$ and let 
$$\Cliff(p^*E) = \Cliff(p^*E)_{\ev} \oplus \Cliff(p^*E)_{\od}$$   denote the complex Clifford algebra bundle \cite{ABS64} of the pullback
$p^*E \to E$ with its natural
$\Z_2$-grading and $G$-action. Consider the smooth $G$-equivariant morphism 
$$C_E : \Cliff(p^*E)_{\ev} \to \Cliff(p^*E)_{\od}$$ 
which is defined via Clifford multiplication by the basepoint
$$C_E(\omega)(e) = e \cdot \omega, \quad e \in E, \omega \in \Cliff(E_{p(e)})_{\ev}.$$   
It is a vector bundle isomorphism off the  zero section $M \subset E$. If $M$ is compact, this defines a $K$-theory element, 
called the Thom class,
$$\lambda_E = [C_E : \Cliff(p^*E)_{\ev} \to \Cliff(p^*E)_{\od}] \in K^0_G(E).$$  
The Thom homomorphism is then defined as the mapping
\begin{eqnarray}\label{eqn:TopThom} 
\Phi : K^0_G(M)  \to  K^0_G(E) : [\sigma]  \mapsto   \lambda_E \otimes [\sigma]
\end{eqnarray} via module multiplication by the Thom class. If
$M$ is not compact, then $\lambda_E$ does not  define a $K$-theory class for $E$, but multiplication by $\lambda_E$ is still well-defined
since $\lambda_E \otimes [\sigma]$ has compact support and we obtain a similar Thom homomorphism in the non-compact case.  An
important result in topological $K$-theory is the Thom isomorphism theorem \cite{Ati67,Kar78,Seg68} which is the statement that if
$p: E \to M$ is an oriented Euclidean vector bundle (of finite even rank) with \spinc-structure then $\Phi: K^0_G(M) \to K^0_G(E)$ is an isomorphism of
abelian groups.

Recall that if we take $M = \{\bullet\}$ to be a point, then $E \cong \R^{2n}$ is a finite dimensional Euclidean vector space and the Thom
isomorphism of $\R^{2n} \to \{\bullet\}$ is the Bott periodicity isomorphism $K^0_G(\{\bullet\}) \cong K^0_G(\R^{2n})$.
The Thom isomorphism can then be viewed, in general, as multiplication by the Bott element $\lambda_{\R^{2n}}$ along the fibers of $E \to M$, where 
$\rank(E) = 2n$.

If a vector bundle $p: \EE \to M$ has {\it infinite} rank, say, if $\EE$ is a real (or complex) Hilbert bundle with infinite dimensional fiber
$\E$, then none of this works. This is because the total space $\EE$ is is not locally compact, and so lies outside the scope of ordinary
topological $K$-theory (even for compact groups.) We could try using an equivariant Swan's Theorem \cite{Ati67,Phi87}
$K^0_G(M) \cong K_0^G(C_0(M))$, where
$C_0(M)$ is the $\cs$-algebra of continuous complex-valued functions on $M$ which vanish at infinity. However, this attempt fails because
$\EE$ is infinite dimensional and so $C_0(\EE) = \{0\}$. Another problem occurs when trying to define an equivariant $K$-theory when the
transformation group $G$ is non-compact, since using finite rank vector $G$-bundles is no longer appropriate. However, N. C. Phillips
\cite{Phi89} has defined a topological equivariant $K$-theory $K^0_G(M)$ for {\it proper} actions of locally compact groups. He then proves a generalized
Green-Julg-Rosenberg isomorphism: $$K^0_G(M) \cong K_0(C_0(M) \rtimes_\alpha G),$$ where $C_0(M) \rtimes_\alpha G$ is (any) crossed product
of $A$ by $G$. But, again we have the problem that $\EE$ is not locally compact.
Most importantly,  non-compact groups do not act {\it properly}  on a point. Thus, we must appeal directly to inherently non-commutative
geometric methods.

In this paper, we remedy these problems by defining an appropriate (equivariant) Thom homomorphism using the $K$-theory of
non-commutative $\cs$-~algebras.  The method is suggested by considering the case when $M = \{\bullet\}$ is a point and so $\EE = \E$ is an infinite
dimensional Euclidean space. In this very important case, we should obtain the Bott periodicity 
of  Higson-Kasparov-Trout \cite{HKT98}. Let $G$ be a smooth, second countable, locally compact Hausdorff topological group, e.g. a countable discrete group.

Let $\E$ be a Euclidean (i.e., real Hilbert) space of countably infinite dimension. If $E_a$ is a finite dimensional subspace of
$\E$, then denote by  $\CC(E_a)$  the $\cs$-algebra $C_0(E_a, \Cliff(E_a))$ of continuous functions from $E_a$ into the complexified
Clifford algebra of $E_a$ vanishing at infinity. This $\cs$-algebra has a $\Z_2$-grading inherited from $\Cliff(E_a)$.  Let $\A(E_a)$
denote the $\Z_2$-graded tensor product,
$$\A(E_a) = C_0(\R) \gtimes \CC(E_a),$$ where $C_0(\R)$ is graded into even and odd functions. If $E_a$ is a subspace of another finite
dimensional subspace $E_b$ of $\E$, then there is a canonical homomorphism of $\cs$-algebras
$\A(E_a) \to \A(E_b)$ (described in more generality in Section 2.) These homomorphisms are injective and,  more
surprisingly, are functorial with respect to subspace inclusion. This allows us to define the direct limit $\cs$-algebra
\begin{eqnarray*}
\A(\E) =\vlim{E_a \subset \E} \A(E_a)
\end{eqnarray*} taken over the directed system of all finite dimensional linear subspaces \{$E_a$\} of the infinite dimensional Euclidean
space $\E$.  If $\E$ is equipped with an orthogonal action of the group $G$, then $G$ acts on the directed system of
subspaces, so it induces an action on the direct limit. Most importantly, the equivariant inclusion of the zero dimensional subspace $\{0\}
\into
\E$ induces a  canonical equivariant $*$-homomorphism $\beta : \A(0) \to \A(\E).$ The Bott periodicity theorem of Higson-Kasparov-Trout
\cite{HKT98} is the statement that the map $\beta : \A(0) \to \A(\E)$ induces an isomorphism in equivariant $K$-theory\footnote{We define $K_*^G(A) = K_*(A \rtimes G)$
where $A \rtimes G$ is the full crossed product (see Definition \ref{def:ktheory}).}: 
$$ \beta_* : K^G_*(\A(0)) \to K^G_*(\A(\E)).$$  Since $\A(0) \rtimes G \cong C_0(\R) \otimes C^*G$ it follows that $K^G_*(\A(0))
\cong K_{*+1}(C^*G)$, where $C^*G$ is the (full) group $\cs$-algebra of $G$. Compare this to the compact group case where
$K^0_G(\{\bullet\}) = R(G) \cong K_0(C^*G)$ by the Peter-Weyl Theorem, where
$R(G)$ is the representation ring of the compact group $G$.

These considerations suggest the following. Let $M$ be a finite dimensional Riemannian manifold. Let $\CC(M)$ denote the $\cs$-algebra of
continuous sections of the complexified Clifford algebra bundle $\Cliff(TM)$ of the tangent bundle $TM$ vanishing at infinity. 
Let $$\A(M) = C_0(\R) \gtimes \CC(M)$$ denote the $\Z_2$-graded suspension of $\CC(M)$.  These non-commutative
$\cs$-algebras will play a fundamental role in defining our Thom homomorphism for
infinite rank Euclidean bundles. If $M$ is an oriented even-dimensional \spinc-manifold, there is a  Morita equivalence \cite{RW98} (given by the spinor bundle) between
$\CC(M)$ and $C_0(M)$. It follows that  if the group $G$ acts properly on $M$ in a compatible way:
$$K_*^G(\A(M)) \cong K_{*+1}(C_0(M) \rtimes G) \cong K^{*+1}_G(M).$$
However, a non-compact group does not act {\it properly} on a point. Thus, we must allow for a more general class of group actions, called {\it isoreductive}
actions (see Definition \ref{def:isoreductive}) which have a local slice property similar to the notion of proper action as given
by Baum-Connes-Higson \cite{BCH94}.

Now let $p : \EE \to M$ be a smooth Euclidean $G$-bundle modelled on $\E$, i.e., a locally trivial fiber bundle with fiber given by a
Euclidean space $\E$, of countably infinite dimension, and structure group $O(\E)$ in the norm (or strong) operator topology. We
assume that the group $G$ acts by Euclidean bundle automorphisms of $\EE$ and acts on $M$ isoreductively.
 We will generalize the $\cs$-algebra construction above to obtain a direct limit
$\cs$-algebra
$$\A(\EE) = \vlim{E^a \subset \EE} \A(E^a)$$ taken over the directed system of all smooth {\it finite rank}
Euclidean subbundles $\EE\supset E^a \to M$ ordered by inclusion of subbundles.  In order for this construction to succeed, we must
equip $\EE$ with a compatible $G$-invariant connection $\nabla$ which controls splittings used in defining the connecting maps of the direct
limit. The equivariant inclusion of $M \into \EE$ as the {\it zero subbundle} then canonically induces an equivariant ``Thom''
$*$-homomorphism  $\Psi_p : \A(M) \to \A(\EE).$ The induced map on equivariant $K$-theory 
$\Psi^G_* : K^G_*(\A(M)) \to K^G_*(\A(\EE))$ is our desired equivariant Thom homomorphism, and our main result is this:
\begin{thm}
Let $G$ be a smooth, second countable, locally compact group and $M$ be a smooth isoreductive Riemannian $G$-manifold. If $p: \EE \to M$ is a
Euclidean $G$-bundle with fiber a countably infinite dimensional Euclidean space $\E$, equipped with a $G$-invariant connection, 
then the map $\Psi_*^G : K^G_*(\A(M)) \to K^G_*(\A(\EE))$ is an isomorphism of abelian groups.
\end{thm}

The outline is as follows. In Section 2, we develop the non-commutative geometric tools needed to define  our infinite rank
Thom homomorphism and discuss the relations with the classical finite rank case. In Section 3, we construct the $\cs$-algebra $\A(\EE)$ 
and prove the non-equivariant version of Theorem 1.1. The equivariant Thom homomorphism is developed in Section 4 and
the proof of Theorem 1.1 is given in Section 5. In Appendix A we review $\Z_2$-graded $\cs$-algebras and their (essentially) self-adjoint
unbounded multipliers and introduce a new balanced tensor product for $\Z_2$-graded $C_0(M)$-algebras that we need to make our construction
work.

The author would like to thank  P. Baum, N. Higson, D. Williams, I. Raeburn, C. Gordon, N. C. Phillips, E. Guentner, and D. Dumitrascu for
helpful suggestions and enlightening discussions.

\section{Clifford $\cs$-algebras and Thom Isomorphism: The Finite Rank Case}

We assume throughout this section that all manifolds are finite dimensional, second countable, Hausdorff,
and infinitely differentiable, and all maps are infinitely smooth. Let $M$ be a Riemannian manifold with tangent bundle $TM \to M$. 

\begin{dfn} 
Let $\Cliff(TM)$ denote the Clifford bundle \cite{ABS64,BGV86} of $TM$. That is, the bundle of Clifford algebras over $M$ whose fiber at
$x \in M$ is the complex Clifford algebra $\Cliff(T_xM)$ of the Euclidean space $T_xM$, i.e., the universal unital complex $*$-algebra that contains
$T_xM$ as a real linear subspace such that $v_x^2 = \|v_x\|^2 1$ for all $v_x \in T_xM$. It has a natural $\Z_2$-graded Hermitian bundle structure.
\end{dfn}  

Recall that the complex Clifford algebra of a finite dimensional Euclidean space has a canonical $\cs$-algebra structure \cite{HKT98,Kar78,LM89}.

\begin{dfn}\label{def:man-alg} 
Denote by $\CC(M)$ the $\cs$-algebra  $C_0(M, \Cliff(TM))$ of continuous sections of $\Cliff(TM)$ which vanish at
infinity, with $\Z_2$-grading induced from $\Cliff(TM)$. $($See Appendix A for a review of graded
$\cs$-algebras.$)$
\end{dfn}

To be precise, we should use the notation $\CC(M, g)$ to show the dependence of the $\cs$-algebra $\CC(M)$ on the Riemannian metric $g$ of
$M$. However, any two bundle metrics $g_1$ and $g_2$ on $TM$ are equivalent via a Euclidean bundle isomorphism $(TM, g_1) \cong (TM, g_2)$ which, by
universality, induces a unitary $\Z_2$-graded bundle isomorphism $\Cliff(TM, g_1) \cong \Cliff(TM, g_2)$ of Clifford algebra bundles.
Thus, the $\cs$-algebras $\CC(M, g_1) \cong \CC(M, g_2)$ are isomorphic as $\Z_2$-graded $\cs$-algebras.

Let $C_0(M)$ denote the $\cs$-algebra of continuous complex-valued functions on $M$ which vanish at infinity. Note that pointwise multiplication
$$(fs)(x) = f(x) s(x), \quad  \forall x \in M,$$ where $f \in C_0(M)$ and $s \in \CC(M)$, determines a nondegenerate $*$-homomorphism
$$C_0(M) \to ZM(\CC(M))$$  into the center of the multiplier algebra of $\CC(M)$ of grading degree zero (where $C_0(M)$ is always considered
to be trivially graded.) Thus, we have the following.

\begin{lem}\label{lem:indep} 
The $\cs$-algebra $\CC(M)$ has a canonical $C_0(M)$-algebra structure, and up to $\Z_2$-graded isomorphism, is
independent of the Riemannian metric on $M$.
\end{lem}

If $M = V$ is a finite dimensional Euclidean vector space, then $TM = V \times V$ and so $\CC(M) = C_0(V, \Cliff(V))$ as in Definition 2.2
\cite{HKT98}. Using the canonical isomorphism
$T(M_1 \times M_2) \cong TM_1 \times TM_2$, we easily obtain the following result.

\begin{lem}\label{lem:trivial} 
If $M = M_1 \times M_2$ then $\CC(M) \cong \CC(M_1) \gtimes \CC(M_2)$.
\end{lem}

\begin{dfn}\label{def:suspend} 
Let $\S$ denote the $\cs$-algebra $C_0(\R)$ of continuous complex-valued functions on the real line which vanish at
infinity, graded into even and odd functions. If $A$ is any $\Z_2$-graded $\cs$-algebra then let $\S A$ be the graded tensor product $\S \gtimes A$.
In particular, let $\A(M) = \S \gtimes \CC(M)$, which can be viewed as a non-commutative topological suspension\footnote{Recall that for the
ungraded tensor product we have $C_0(\R) \otimes C_0(M) \cong C_0(\R \times M).$} of $M$.
\end{dfn}

Another relationship between the non-commutative $\cs$-algebra $\CC(M)$ and the commutative $\cs$-algebra
$C_0(M)$ is given by \spinc-structures \cite{LM89}. Let $\C_1 = \Cliff(\R)$ denote the first complex Clifford algebra. The following is adapted from
Theorem 2.11 of Plymen \cite{Ply86} and Proposition II.A.9 of Connes \cite{Con94}.

\begin{prop}\label{prop:morita} 
Let $M$ be an oriented Riemannian manifold. If $M$ is even-dimensional, there is a bijective correspondence between
\spinc-structures on $M$ and Morita equivalences $($in the sense of Rieffel \cite{RW98,Rie82}$)$ between the $\cs$-algebras $C_0(M)$ and
$\CC(M)$. Thus, $\A(M)$ is Morita equivalent to $C_0(\R \times M)$. If $M$ is odd-dimensional, then
\spinc-structures on $M$ are in bijective correspondence with Morita equivalences $C_0(M) \sim \CC(M) \gtimes \C_1$.
\end{prop}

\Pf Assume $\dim(M) = 2m$ is even. Let $S \to M$ denote the spinor bundle associated to a \spinc-structure. This has a
natural Euclidean bundle structure. The module $\A_S = C_0(M,S)$ of continuous sections of $S$ vanishing at infinity gives the desired
$C_0(M)-\CC(M)$-imprimitivity bimodule implementing the Morita equivalence $C_0(M) \sim \CC(M)$. Since $M$ is even-dimensional, 
$$\Cliff(TM) = \Cliff(TM)_{\ev} \oplus \Cliff(TM)_{\od}$$ splits as a $\Z_2$-bundle with grading operator given by the {\it bounded} global section
$$\epsilon = i^me_1 \cdot e_2 \cdots e_{2m},$$ where $\{e_1, \cdots, e_{2m}\}$ is any oriented local orthonormal frame of $TM$. It follows that
$\CC(M) = C_0(M, \Cliff(TM))$ is an {\it evenly} graded $\cs$-algebra with grading operator $\epsilon \in M(\CC(M))$. Thus, by Proposition 14.5.1
\cite{Bla98}, the  graded tensor product $C_0(\R) \gtimes \CC(M)$ and ungraded tensor product $C_0(\R) \otimes \CC(M)$ are isomorphic as (ungraded)
$\cs$-algebras. Thus, since suspensions preserve Morita equivalence:
$$\A(M) =_{\text{def}} C_0(\R) \gtimes \CC(M) \cong C_0(\R) \otimes \CC(M) \sim C_0(\R) \otimes C_0(M) \cong C_0(\R \times M).$$ 
The odd-dimensional case is similar, but we cannot remove the grading on the suspension and we must tensor by $\C_1$.
\EPf

Although $\CC(M)$ and $\A(M)$ carry natural $\Z_2$-gradings, when we consider their $K$-theory, we shall ignore the gradings. To make
this point clearer, {\bf if $A$ is any $\cs$-algebra---graded or not---then $K_*(A)$ will denote the $K$-theory of the underlying
$\cs$-algebra, without the grading.}

\begin{cor}\label{cor:topKthry} 
If $M$ is an even-dimensional oriented Riemannian manifold with \spinc-structure, there are canonical $K$-theory isomorphisms
$$K_*(\CC(M)) \cong K^*(M) \text{ and }K_*(\A(M)) \cong K^{*+1}(M).$$
\end{cor}

Let $p : E \to M$ be a smooth finite rank Euclidean bundle. We wish to show that there is
a natural ``Thom'' $*$-homomorphism
$$\Psi_p : \A(M) \to \A(E)$$  where we consider $E$ as a finite dimensional manifold with Riemannian structure to be constructed as follows.

Given $p : E \to M$, there is a short exact sequence \cite{AMR88,BGV86} of real vector bundles 
$$\xymatrix{ 0 \ar[r] & VE \ar[r] & TE \ar[r]^-{T^*p} & p^*TM \ar[r] & 0}$$  
where the {\it vertical subbundle} $VE = \ker(T^*p)$ is isomorphic to $p^*E$. This
sequence does {\it not} have a canonical splitting, in general, but choosing a compatible connection $\nabla$ on $E$ determines an associated vector bundle
splitting.

\begin{dfn} 
A {\bf connection} \cite{BGV86,LM89} on $E$ is a linear map
$$\nabla : C^\infty(M, E) \to C^\infty(M, T^*M \otimes E)$$  which satisfies Leibnitz's rule, i.e., if $s \in C^\infty(M,E)$ and $f \in C^\infty(M)$,
then
$$\nabla(fs) = df \otimes s + f\nabla s$$  where $d$ is the exterior derivative on smooth forms $\Omega^*(M) = C^\infty(M, \Lambda^*T^*M)$. If $E$ has
a Euclidean metric $(\cdot ,\cdot )$, we say that
$\nabla$ is {\bf compatible} with the metric if 
$$X(s_1,s_2) = (\nabla_X s_1, s_2) + (s_1, \nabla_X s_2)$$  for all sections $s_1, s_2 \in C^\infty(M,E)$ and vector fields $X \in C^\infty(M, TM)$.
If $E$ is equipped with a compatible connection $\nabla$, then we call $E$ an {\bf affine} Euclidean bundle.
\end{dfn}

\begin{lem}\label{lem:splitting} 
Let $p : E \to M$ be a finite rank Euclidean bundle on the Riemannian manifold $M$. Given a compatible connection
$\nabla$ on $E$ there is an induced orthogonal splitting $TE \cong p^*E \oplus p^*TM$ of the exact sequence
$$0 \to p^*E \to TE \to p^*TM \to 0$$ where $p^*E$ and $p^*TM$ have the pullback metrics. Hence, the manifold $E$ has an induced Riemannian metric.
\end{lem}

\Pf Let $\nabla^{*} :C^\infty(E, p^*E) \to C^\infty(T^*E \otimes p^*E)$ denote the pullback of $\nabla$ on the bundle $p^*E \to E$, which is defined
by the formula:
$$\nabla^{*}(f p^*s) = df \otimes p^*s + f p^*(\nabla s)$$ for $f \in C^\infty(M)$ and $s \in C^\infty(M,E)$. The tautological section $\tau \in
C^\infty(E, p^*E)$ is the smooth section defined by the formula $\tau(e) = (e, e)$ for all $e \in E$. The derivative of $\tau$ will be denoted by
$$\omega = \nabla^{*} \tau \in C^\infty(T^*E \otimes p^*E) = \Omega^1(E, p^*E) \cong \Omega^1(E, VE)$$ 
which is a connection $1$-form (Definition 1.10 \cite{BGV86}). The kernel $HE = \ker(\omega) \cong p^*TM$ of the connection
$1$-form $\omega$ is the horizontal subbundle of $TE$ which provides the desired splitting
$$TE = VE \oplus HE \cong p^*E \oplus p^*TM.$$  Now give $TE$ the direct sum of the pullback metrics on $p^*E$ and $p^*TM$. This gives $E$ the
structure of a Riemannian manifold and makes the splitting of $TE$ orthogonal. 
\EPf

Thus, given a compatible connection $\nabla$ on the Euclidean bundle $E$, we can define the $\cs$-algebra $\CC(E)$ as above using the induced
Riemannian structure on the {\it manifold} $E$. However, we also have the $\cs$-algebra $C_0(E, \Cliff(p^*E))$ associated to the pullback bundle
$p^*E \to E$. Both $\CC(E)$ and $C_0(E, \Cliff(p^*E))$ have natural
$C_0(E)$-algebra structures. However, the bundle map $p: E \to M$ induces a pullback
$*$-homomorphism \cite{RW85} $$p^* : C_0(M) \to C_b(E) = M(C_0(E)) : f \mapsto p^*(f) = f \circ p$$  which induces a (graded) $C_0(M)$-algebra structure on any
(graded) $C_0(E)$-algebra. The following is an important result that relates these two $\cs$-algebras\footnote{Although $C_0(E, \Cliff(p^*E)) \cong p^*C_0(M, \Cliff(E))$,
we will not need this isomorphism.} to the $C_0(M)$-algebra $\CC(M)$ of the base manifold $M$.

\begin{thm}\label{thm:decomp} 
Let $p : E \to M$ be a finite rank affine Euclidean bundle on the Riemannian manifold $M$.
There is a natural isomorphism of graded
$\cs$-algebras
$$\CC(E) \cong C_0(E, \Cliff(p^*E)) \gtimes_{_M} \CC(M)$$
where $\gtimes_{_M}$ denotes the balanced tensor product over $C_0(M)$ of Definition \ref{def:gtimes_M}. 
\end{thm}

\Pf By the previous lemma, there is an induced orthogonal splitting 
$$TE = p^*E \oplus p^*TM.$$   Thus,  we have an induced isomorphism of $\Z_2$-graded Clifford algebra bundles
\begin{eqnarray}\label{eqn:cliffbdl}
\Cliff(TE) \cong \Cliff(p^*E \oplus p^*TM) = \Cliff(p^*E) \gtimes p^*\Cliff(TM).
\end{eqnarray} Therefore, by taking sections, we have canonical balanced tensor product isomorphisms (see Proposition \ref{prop:balanced})
\begin{eqnarray*}
\CC(E) & =_{\text{def}}  & C_0(E, \Cliff(TE)) \cong C_0(E, \Cliff(p^*E) \gtimes p^*\Cliff(TM)) \\
       & \cong_{\phantom{\text{def}}} & C_0(E, \Cliff(p^*E)) \gtimes_{_E} C_0(E, \Cliff(p^*TM)).
\end{eqnarray*} But, we have, using pullbacks along $p : E \to M$, that there are canonical pullback isomorphisms 
(Proposition \ref{prop:pullbackalg})
$$C_0(E, \Cliff(p^*TM)) \cong p^*C_0(M, \Cliff(TM)) = p^*\CC(M) =_{\text{def}} C_0(E) \gtimes_{_M} \CC(M).$$ 
Hence, it follows that
\begin{eqnarray*}
\CC(E) & \cong_{\phantom{\text{def}}} & C_0(E, \Cliff(p^*E)) \gtimes_{_E} C_0(E, \Cliff(p^*TM)) \\
       & \cong_{\phantom{\text{def}}} & C_0(E, \Cliff(p^*E)) \gtimes_{_E} C_0(E) \gtimes_{_M} \CC(M) \\
       &  \cong_{\phantom{\text{def}}} & C_0(E, \Cliff(p^*E)) \gtimes_{_M} \CC(M)
\end{eqnarray*} using the canonical isomorphism $A \gtimes_{_E} C_0(E) \cong A$ for graded $C_0(E)$-algebras. \EPf 

We now wish to define a certain ``Thom operator'' for the ``vertical'' algebra $C_0(E, \Cliff(p^*E))$.

Associate to the Euclidean bundle $E$ an unbounded section
\begin{eqnarray*} C_E : E  \to  \Cliff(p^*E) : e & \mapsto C_{p(e)}(e) \end{eqnarray*}   where $C_{p(e)}$ is the Clifford operator on the Euclidean
space $E_{p(e)}$ from Definition 2.4 \cite{HKT98}. It is given globally by the composition
$$\xymatrix{ E \ar[r]^-{\tau} \ar@/_1.5pc/[rr]_-{C_E}  & p^*E \ar[r]^-{C} & \Cliff(p^*E)}$$   where $\tau \in C^\infty(E, p^*E)$ is the tautological section
(in the proof of Lemma \ref{lem:splitting}) and
$C : p^*E \into \Cliff(p^*E)$ is the canonical inclusion $C(e_1, e_2) = C_{p(e_1)}(e_2)$.   The following is then easy to prove.

\begin{thm}\label{thm:cliffop} 
Let $E$ be a finite rank Euclidean bundle on $M$. Multiplication by the section $C_E : E \to \Cliff(p^*E)$  determines a
degree one, essentially self-adjoint, unbounded multiplier $($see Definition \ref{def:usam}$)$ of the $\cs$-algebra $C_0(E, \Cliff(p^*E))$ with domain $C_c(E,
\Cliff(p^*E))$.
\end{thm}

We will call $C_E$ the {\bf Thom operator} of $E \to M$. Thus, we have a functional calculus homomorphism
$$\S = C_0(\R) \to M(C_0(E, \Cliff(p^*E))) : f \to f(C_E)$$  from $\S$ to the multiplier algebra of $C_0(E, \Cliff(p^*E))$. Note that $f(C_E)$ goes to
zero in the ``fiber'' directions on $E$ (since $p(e)$ is constant), but is only bounded in the ``manifold'' directions on $E$.  Indeed, note that for
the generators $f(x) = \exp(-x^2)$ and $g(x) = x \exp(-x^2)$ of $\S$, we have that $f(C_E)$ and $g(C_E)$ are, respectively, multiplication by the
following functions on $E$:
$$f(C_E)(e) = \exp(-\|e\|^2) \text{ and } g(C_E)(e) = e \cdot \exp(-\|e\|^2), \quad \forall e \in E.$$ 

\begin{dfn} Let $X$ denote the degree one, essentially self-adjoint, unbounded multiplier of $\S = C_0(\R)$, with domain the compactly supported
functions, given by multiplication by $x$, i.e., $Xf(x) = xf(x)$ for all $f \in C_c(\R)$ and $x \in \R$. \end{dfn}

By Lemma \ref{lem:tensormult}, the operator $X \gtimes 1 + 1 \gtimes C_E$ determines a degree one, essentially self-adjoint, unbounded multiplier of 
the tensor product
$$\S \gtimes C_0(E, \Cliff(p^*E)) = \S C_0(E, \Cliff(p^*E)).$$  We obtain a functional calculus homomorphism 
$$\beta_E : \S \to M(\S C_0(E, \Cliff(p^*E))) : f \mapsto f(X \gtimes 1 + 1 \gtimes C_E)$$
from $\S$ into the multiplier algebra of $\S C_0(E, \Cliff(p^*E))$.  Now we can define our ``Thom $*$-homomorphism'' for a
finite rank affine Euclidean bundle. This will provide the connecting map in the next section to define the direct limit $\cs$-algebra for an infinite rank
affine Euclidean bundle.

\begin{thm}\label{thm:Thomhom} 
Let $(E,\nabla,M)$ be as above. With respect to the isomorphism
$$\A(E) \cong \S \gtimes C_0(E, \Cliff(p^*E)) \gtimes_{_M} \CC(M)$$ from Theorem \ref{thm:decomp}, there is a graded $*$-homomorphism 
$$\Psi_p  =  \beta_E \gtimes_{_M} \id_M : \A(M) \to \A(E)$$    which on elementary tensors $f \gtimes s \in \S \gtimes \CC(M)=\A(M)$ is given by
$$f \gtimes s \mapsto f(X \gtimes 1 + 1 \gtimes C_E) \gtimes_{_M} s.$$
\end{thm}

\Pf From the discussion above, we have that $\beta_E \gtimes_{_M} \id_M$ is the composition
$$\xymatrix{
\S \gtimes \CC(M) \ar[r]^-{\beta_E \gtimes \id_M} & M(\S C_0(E, \Cliff(p^*E))) \gtimes \CC(M) \to M(\S C_0(E, \Cliff(p^*E))) \gtimes_{_M} \CC(M)
}$$
Checking on the generator $f(x) = \exp(-x^2)$ of $\S$, we compute that
$$f(X \gtimes 1 + 1 \gtimes C_E) \gtimes_{_M} s = \exp(-x^2) \gtimes \exp(-\|e\|^2) \gtimes_{_M} s \in \A(E)$$ Similarly for $g(x) = x \exp(-x^2)$,
we find that
\begin{eqnarray*} g(X \gtimes 1 + 1 \gtimes C_E) \gtimes_{_M} s & =  &x\exp(-x^2) \gtimes \exp(-\|e\|^2) \gtimes_{_M} s \\
 & + & \exp(-x^2) \gtimes e \cdot \exp(-\|e\|^2)\gtimes_{_M} s\in \A(E).
\end{eqnarray*}  It follows that the range of $\Psi_p = \beta_E \gtimes_{_M} \id_M$ is in $\A(E)$ as desired.
\EPf

\noindent {\bf Important Note:} Compare the similarities of the $\cs$-algebraic formula above to formula (\ref{eqn:TopThom}) for the topological Thom
homomorphism $\Phi : K^0(M) \to K^0(E)$. Both have the same form: pullback from $M$ and multiply/tensor by the Thom class/operator of the bundle $E$.

We now come to the main theorem of this section:

\begin{ThomThm}\label{thm:Thomfrc} 
If $E \to M$ is a smooth finite rank affine Euclidean bundle, then the $*$-homomorphism $\Psi_p : \A(M) \to \A(E)$ induces an isomorphism  of abelian groups:
$$\Psi_{*} : K_*(\A(M)) \to K_*(\A(E))$$
\end{ThomThm}

The proof will be essentially the same argument as given in another ``Thom isomorphism'' theorem of the author (Theorem B.22 \cite{Trou99}). But first,
we must develop some functorial properties of these Thom maps.

Since the space of compatible connections $\nabla$ on $E \to M$ is convex, we have the following result.

\begin{prop}\label{prop:htpy} 
Let $p : E \to M$ be a smooth finite rank affine Euclidean bundle. The homotopy class of the  $*$-homomorphism $\Psi_p : \A(M) \to
\A(E)$ is independent of the choice of compatible connection $\nabla$ on $E$. 
\end{prop}

\begin{prop}\label{prop:triv} 
If $p : E = M \times V \to M$ is a trivial finite rank affine Euclidean bundle $($with trivial connection $\nabla_0 = d$ $)$ then we
have a $\Z_2$-graded isomorphism
$$\CC(E) \cong \CC(V) \gtimes \CC(M)$$ such that the Thom map has the form
$$\Psi_p \cong \beta_V \gtimes \id_{\CC(M)} : \A(M) = \S \gtimes \CC(M) \to \A(V) \gtimes \CC(M) \cong \A(E)$$   where $\beta_V : \S \to \A(V) : f \mapsto
f(X \gtimes 1 + 1 \gtimes C_V)$ is the Thom map for $V \to \{0\}$.
\end{prop}

\Pf The trivial connection $\nabla_0 = d$ gives the manifold $E = M \times V$ the Riemannian metric induced by the isomorphism 
$$TE = TM \times TV \to M \times V =  E.$$ The pullback vector bundle $p^*E \to E$ has the form
$$p^*E = (M \times V) \times V \to M \times V = E$$ and so the Clifford bundle $\Cliff(p^*E) = (M \times V) \times \Cliff(V)$, which gives:
$$C_0(E, \Cliff(p^*E)) = C_0(M \times V, (M \times V) \times \Cliff(V)) \cong C_0(V, \Cliff(V)) \gtimes C_0(M).$$ 
By Theorem \ref{thm:decomp}, it
follows that
\begin{eqnarray*}
\CC(E) & \cong & C_0(E, \Cliff(p^*E)) \gtimes_{_M} \CC(M) \\
       & \cong & C_0(V, \Cliff(V)) \gtimes C_0(M) \gtimes_{_M} \CC(M) \\
       & \cong & \CC(V) \gtimes \CC(M).
\end{eqnarray*} 
where we used the isomorphism $C_0(M) \gtimes_{_M} \CC(M) \cong \CC(M)$. The result now easily follows. \EPf

For example, if $p: E_b \to E_a$ is the orthogonal projection of a finite dimensional Euclidean vector space $E_b$ onto a linear subspace $E_a$ then
$\Psi_p = \beta_{ba}$ is the  ``Bott homomorphism'' from Definition 3.1 of Higson-Kasparov-Trout \cite{HKT98}. 

\begin{lem}\label{lem:open} 
Let $U$ be an open subset of the Riemannian manifold $M$. The inclusion $i : U \into M$ induces a short exact sequence
$$\xymatrix{ 0 \ar[r] & \A(U) \ar[r]^-{1 \gtimes i_*}  & \A(M) \ar[r] &  \A(M \minus U) \ar[r] & 0
}$$ of $\cs$-algebras. Thus, $\A(U) \ideal \A(M)$ as a $($two-sided$)$ $\cs$-ideal.
\end{lem}

\Pf Give the open subset $U$ the restriction $g_U = i^*(g)$ of the Riemannian metric $g$ on $M$. The tangent bundle
$TU = TM|_U = i^*(TM)$ is the restriction of $TM$ to $U$. Thus, it follows by universality that
$$\Cliff(TU) = \Cliff(i^*(TM)) = i^*\Cliff(TM) = \Cliff(TM)|_U.$$ Since $\CC(U) = C_0(U, \Cliff(TU))$ we have an exact sequence of $\cs$-algebras
$$\xymatrix{ 0 \ar[r] & \CC(U) \ar[r]^-{i_*}  & \CC(M) \ar[r] &  \CC(M \minus U) \ar[r] & 0.}$$ 
Taking the (maximal) graded tensor product with the
nuclear, hence exact \cite{Wass94}, $\cs$-algebra $\S = C_0(\R)$, we obtain an exact sequence of $\cs$-algebras
$$\xymatrix{ 0 \ar[r] & \S \gtimes C(U) \ar[r]^-{1 \gtimes i_*}  & \S \gtimes C(M) \ar[r] &  \S \gtimes C(M \minus U) \ar[r] & 0
}$$ as desired. \EPf

The next result shows that the functor $M \mapsto K_*(\A(M))$ has Mayer-Vietoris exact sequences.

\begin{lem}\label{lem:mayer} 
Let $U$ and $V$ be open subsets of $M$. There is an exact sequence
$$\xymatrix{
K_0(\A(U \cap V)) \ar[r] & K_0(\A(U)) \oplus K_0(\A(V)) \ar[r] & K_0(\A(U \cup V)) \ar[d] \\
K_1(\A(U \cup V)) \ar[u] & K_1(\A(U)) \oplus K_1(\A(V)) \ar[l] & K_1(\A(U \cap V)) \ar[l]
}$$
of abelian groups. 
\end{lem}

\Pf By Lemma \ref{lem:open}, $\A(U)$ and $\A(V)$ are ideals in $\A(U \cup V)$. By an approximation argument we have the
relations:
\begin{eqnarray*} \Bigg\{
\aligned
\A(U) + \A(V) & = \A(U \cup V) \\
\A(U) \cap \A(V) & = \A(U \cap V)
\endaligned
\end{eqnarray*}
since similar relations hold for $\CC(U)$ and $\CC(V)$. The result now follows from Exercise 4.10.21 \cite{HR00}.
\EPf

\begin{lem}\label{lem:openfunc} 
Let $(E, \nabla, M)$ be as above. If $U$ is an open subset of $M$ there is a commutative  diagram
$$\CD
\A(E|_U) @>{1 \gtimes \tilde{i_*}}>> \A(E) \\ @A{\Psi_{p_{_U}}}AA @AA{\Psi_p}A  \\
\A(U) @>{1 \gtimes i_*}>> \A(M)
\endCD$$ where $p_U : E|_U \to U$ denotes the restriction of the Euclidean bundle $E$ to $U$.
\end{lem}

\Pf Let $i : U \into M$ denote the inclusion of $U$ as an open submanifold of $M$. The restricted bundle $E|_U = i^*(E)$, which as a subset of $E$, is
open since $E|_U = p^{-1}(U)$ as a set. We have a commuting diagram
$$\CD  E|_U @>{\tilde{i}}>> E \\ @V{p_U}VV @VV{p}V \\ U @>{i}>> M 
\endCD$$ where the horizontal maps are the inclusions as open subsets and the vertical maps are the bundle projections. Give $E|_U$ the restriction of
the Euclidean metric on $E$. Let $\nabla_U = i^*(\nabla)$ denote the restriction of $\nabla$ to $E|_U$.  The tautological section $\tau_U : E|_U \to
p^*E|_U$ is the restriction $\tau_U = \tau|_{E_U} = \tilde{i}^*(\tau)$ of the tautological section $\tau$ of $E$. It follows that the connection
$1$-form $\theta_U = \nabla^*_U \tau_U \in \Omega^1(E|_U, p^*E|_U)$ satisfies
$$\theta_U = \nabla^*_U \tau_U = p_U^*i^*(\nabla) \tilde{i}^*(\tau) =\tilde{i^*}(\nabla \tau) = \tilde{i^*}(\theta)$$ and so is the restriction of the
connection $1$-form $\theta \in \Omega^1(E, p^*E)$ of $\nabla$ to $E|_U$. This implies that the horizontal subbundle $H(E|_U)$ is
$$p_U^*TU \cong H(E|_U) = \ker{\theta_U} \cong \tilde{i}^*(\ker{\theta}) = \tilde{i}^*(HE) \cong (p^*TM)|_{E|_U}.$$ Thus, we have coherent splittings
$$\xymatrix{ TE  & = & p^*E \oplus TM \\ TE|_U \ar[u]^{T\tilde{i}} &  = & p^*E|_U \oplus TU \ar[u]_{p^*\tilde{i} \oplus Ti}
}$$ where the vertical maps are open inclusions. This implies that we have an induced commuting diagram
$$\xymatrix{
\CC(E)  & \cong & C_0(E, \Cliff(p^*E)) \gtimes_M \CC(M) \\
\CC(E|_U) \ar[u]^-{\tilde{i}_*} & \cong & C_0(E|_U, \Cliff(p^*E|_U)) \gtimes_U \CC(U) \ar[u]_-{\tilde{i}_* \gtimes i_*}
}$$ of $\cs$-algebras. Since the Thom operator $C_{E_U} : E|_U \to \Cliff(p^*E|_U)$ is the restriction of the Thom operator $C_E$
to $E|_U$,the result easily follows by checking on generators $\exp(-x^2) \gtimes s_U$ and $x \exp(-x^2) \gtimes s_U$ of $\S \gtimes \CC(U)$.
\EPf

We are now ready to prove the finite rank Thom Isomorphism Theorem \ref{thm:Thomfrc}.

\smallskip

\Pff {\it of Theorem \ref{thm:Thomfrc}}. If $M = \{0\}$ is a one-point space, then $E = V$ is a finite dimensional Euclidean space and so, by Lemma
\ref{lem:trivial}, $\Psi_p = \beta_V : \S = \A({0}) \to \A(V).$ The induced map on $K$-theory $\Psi_* = \beta_* : K_*(\S) \to K_*(\A(V))$ is an
isomorphism by the Bott periodicity theorem 2.6 of Higson-Kasparov-Trout \cite{HKT98}. 

If $E = M \times V \to M$ is a trivial bundle, then we may assume that $\nabla = d$ is the trivial connection by Proposition
\ref{prop:htpy}. Thus, by Lemma \ref{lem:trivial} we have that
$\Psi_p = \beta_V \gtimes \id_{\CC(M)} : \A(M) = \S \gtimes \CC(M) \to \A(V) \gtimes \CC(M) \cong \A(E)$. The induced map on $K$-theory is an isomorphism, again by the
proof of Theorem 2.6 \cite{HKT98}, by tensoring with $\id_{\CC(M)}$. Indeed, the inverse of $\Psi_*$ is induced by
$$\{\phi_t \gtimes \id_{\CC(M)}\} : \A(V) \gtimes \CC(M) \to \S \K(L^2(V, \Cliff(V)) \gtimes \CC(M)$$
where 
$$\{\phi_t\} : \A(V) \to \S \K(L^2(V, \Cliff(V)) : f \gtimes s \mapsto f(X \gtimes 1 + 1 \gtimes t^{-1}D_V)(1 \gtimes s)$$
is the asymptotic morphism of the Dirac operator $D_V$ on $L^2(V, \Cliff(V))$ of Definitions 2.7 \& 2.8 \cite{HKT98}.

Let $\{U_n\}_1^\infty$ be a countable covering by open subsets of $M$ such that $E$ trivializes over each $U_n$. We now proceed by induction on $n$.  If $n=1$ the
result follows by Lemma \ref{lem:indep}, Proposition \ref{prop:triv} and the discussion above. Suppose we have proved the result for $E$ restricted to
$U = U_1 \cup U_2 \cup \cdots \cup U_n$ (or any open subset of $U$). Set $V = U_{n+1}$. By Lemmas
\ref{lem:open} and \ref{lem:openfunc}, there is a commuting cubical diagram of
$\cs$-algebra homomorphisms
$$\xymatrix{
\A(E|_U) \ar[rr]   &  & \A(E|_{U \cup V})   & \\
 & \A(E|_{U \cap V}) \ar[ul] \ar[rr]  & & \A(E|_V) \ar[ul]  \\
\A(U) \ar[rr] \ar[uu] & & \A(U \cup V) \ar[uu]  & \\
 & \A(U \cap V) \ar[ul] \ar[uu] \ar[rr] & & \A(V)  \ar[ul] \ar[uu] 
}$$  where the vertical arrows are the appropriate Thom maps and the top and bottom squares are inclusions. By Lemma \ref{lem:mayer} 
we have Mayer-Vietoris sequences
$$\xymatrix@C-1pc{
\cdots \ar[r] & k^*(E|_{U \cap V}) \ar[r]  & k^*(E|_U)) \oplus k^*(E|_V) \ar[r] & k^*(E|_{U \cup V}) \ar[r] & k^{*+1}(E|_{U \cap V})  \ar[r] & \cdots \\
\cdots \ar[r] & k^*(U \cap V) \ar[r] \ar[u] & k^*(U) \oplus k^*(V) \ar[r] \ar[u] & k^*(U \cup V) \ar[r] \ar[u] & k^{*+1}(U \cap V) \ar[u] \ar[r] & \cdots \\
}$$ 
where we have used the simplifying notation $k^*(M) = K_*(\A(M))$. By induction, the vertical arrows are isomorphisms from $k^*(U)$,
$k^*(V)$, and $k^{*}(U \cap V)$ (since $U \cap V$ is an open subset of $U$.) By the Five Lemma, we conclude that the vertical arrows from $k^{*}(U \cup V)$ 
are also isomorphisms. This concludes the proof.
\EPf

\begin{cor} 
If $E$ is a finite even rank oriented Euclidean \spinc-bundle $($with spin connection $\nabla$$)$ on an even-dimensional
oriented Riemannian \spinc-manifold $M$, then $\Psi_p : \A(M) \to \A(E)$ induces the topological Thom isomorphism 
$$\xymatrix{
 K_*(\A(M)) \ar[r]^{\Psi_*} \ar[d]_{\cong} & K_*(\A(E)) \ar[d]^{\cong} \\
 K^{*+1}(M) \ar[r]^{\Phi} & K^{*+1}(E)
}$$
\end{cor}

\Pf From the splitting $TE = p^*E \oplus p^*TM$, we see that $E$ has a canonical structure as an oriented even-dimensional \spinc-manifold. The result now follows by
Proposition \ref{prop:morita}, Corollary \ref{cor:topKthry}, and the previous proof. \EPf

\section{The $\cs$-algebra of an Infinite Rank Euclidean Bundle}

In this section, we are going to construct a $\cs$-algebra $\A(\EE)$ for an infinite rank Euclidean bundle $\EE \to M$, analogous to the one in
Definition \ref{def:suspend}. It should be stressed that in the infinite dimensional case $\A(\EE)$ is {\it not} the suspension of a $\cs$-algebra. 
We choose the notation $\A(\EE)$ to be consistent with the notation in \cite{HK01}, which differs from the potentially confusing notation
in the original Bott periodicity paper \cite{HKT98}.
First we need to investigate the transitivity properties of the Thom maps from the previous section.

Fix a smooth finite dimensional Riemannian manifold $M$. Let $p_2 : E^2 \to M$ be a finite rank Euclidean bundle with
metric $( \cdot, \ cdot)_2$. Suppose that $E^1 \subset E^2$ is a vector subbundle of $E^2$, i.e., there is an exact sequence of vector bundles
$$\xymatrix{
0 \ar[r] & E^1 \ar[r]^-{i} & E^2.
}$$
Let $p_1 = p_2|_{E_1} : E^1 \to M$ be the subbundle projection. 

We can restrict the bundle metric on $E^2$ to get a metric $( \cdot , \cdot)_1 = i^*( \cdot, \cdot )_2$ on $E^1$ so that $i$ is an inclusion of Euclidean bundles. We
then have a canonical orthogonal splitting
$E^2 = E^\perp \oplus E^1$ of Euclidean bundles. Now let $\nabla^2$ be a compatible connection on $E^2$.
We can then define a connection $\nabla^1$ on $E^1$ by the composition
$$\xymatrix{  C^\infty(M, E^2) \ar[r]^-{\nabla^2} & C^\infty(M, T^*M \otimes E^2) \ar[d]^-{(1 \otimes p_{21})_*} \\
 C^\infty(M, E^1) \ar[u]^-{i_*} \ar[r]^-{\nabla^1} & C^\infty(M, T^*M \otimes E^1)
}$$ 
where $p_{21} : E^2 = E^\perp \oplus E^1 \to E^1$ is the orthogonal bundle projection onto $E^1$. The connection $\nabla^1$ is the called the {\it
projection} of the connection $\nabla^2$ onto $E^1$.

\begin{lem} 
The projected connection $\nabla^1$ on the subbundle $E^1$ is compatible with the induced metric $( \cdot , \cdot)_1$.
\end{lem}

Let $E^{21}$ denote the vector bundle $p_{21} : E^2 \to E^1$. There is then a Euclidean bundle isomorphism $E^{21} \cong p_1^*E^\perp$, and  an induced
connection $\nabla^{21} = p_1^*\nabla^\perp$, where $\nabla^\perp$ is the induced connection on $p_2 : E^\perp \to M$ via projecting
$\nabla^2$ as above. (Note that $\nabla^2 = \nabla^\perp \oplus \nabla^1$ with respect to the orthogonal splitting $E^2 = E^\perp \oplus E^1$.)
Hence, we have a commuting diagram of vector bundle projections
$$\xymatrix{ & E^2 \ar[dd]^{p_2} \ar[dl]_{p_{21}} \\
   E^1 \ar[dr]_{p_1} & \\ & M }$$ 
and an associated diagram of Thom $*$-homomorphisms (via Theorem \ref{thm:Thomhom}):
\begin{eqnarray}
\xymatrix{
 & \A(E^2) \\
\A(E^1) \ar[ur]^{\Psi_{21}} &  \\ & \A(M) \ar[ul]^{\Psi_1} \ar[uu]_{\Psi_2}
}
\end{eqnarray} 
We want to prove that this diagram commutes, i.e., the Thom maps are transitive. The proof is similar to Proposition 3.2 \cite{HKT98}, which will be an
immediate corollary.

\begin{prop}\label{prop:transitive} 
Let $p_2 : E^2 \to M$ be a finite rank Euclidean bundle with connection $\nabla^2$. If $p_1 : E^1 \to M$ is a
finite rank vector bundle such that there is a vector bundle inclusion $i: E^1 \hookrightarrow E^2$, then diagram (1) above commutes, where $\Psi_{21} :
\A(E^1) \to \A(E^2)$  is the Thom homomorphism associated to the orthogonal projection $p_{21} : E^2 \to E^1$.
\end{prop}
 
\Pf First we need to check that the Riemannian metric (from Lemma \ref{lem:splitting}) on the total space manifold $E^2$ induced by the triple 
$(E^2, \nabla^2, M)$  (used to define the $\cs$-algebra $\CC(E^2)$) is the same as the Riemannian metric on the manifold $E^{21} = E^2$ induced by the
triple  $(E^{21}, \nabla^{12}, E^1)$. Using the induced metrics and connections from above, we have the following orthogonal splittings of bundles:
$$\left\{ \aligned 
TE^2 & =  p_2^*E^2 \oplus p_2^*(TM) \\
TE^{21} & =  p_{21}^*E^{21} \oplus p_{21}^*(TE^1) \\
TE^1&  =   p_1^*E^1  \oplus p_1^*(TM)
\endaligned \right.$$
By plugging the third into the second and using the fact that $p_2 = p_1 \circ p_{21}$, we obtain the first splitting. Hence, these
orthogonal splittings are compatible with each other. Therefore, there is a canonical identification $\CC(E^2) = \CC(E^{21})$ of graded $\cs$-algebras.

The orthogonal splitting $E^2 = E^\perp \oplus E^1$ induces an orthogonal splitting
$$p_2^*E^2 = p_2^*E^\perp \oplus p_2^*E^1 = p_{21}^* E^{21} \oplus p_{21}^* p_1^*E^1.$$
This implies that we have a canonical Clifford algebra bundle decomposition
\begin{eqnarray}\label{eqn:cliffbdl2}
\Cliff(p_2^*E^2) \cong \Cliff(p_{21}^*E^{21}) \gtimes p_{21}^*\Cliff(p_1^*E^1).
\end{eqnarray}
This decomposition induces a canonical isomorphism of graded $\cs$-algebras
\begin{eqnarray*}
C_0(E^2, \Cliff(p_2^*E^2)) & \cong & C_0(E^2, \Cliff(p_{21}^*E^{21})) \gtimes_{_{E^2}} p_{21}^*C_0(E^1, \Cliff(p_1^*E^1)) \\
& = & C_0(E^{21}, \Cliff(p_{21}^*E^{21})) \gtimes_{_{E^2}}   p_{21}^*C_0(E^1, \Cliff(p_1^*E^1))
\end{eqnarray*}
where we use the fact that $E^{21} = E^2$ as {\it manifolds}. Consequently, we have
\begin{eqnarray}\label{eqn:c*decomp}
\A(E^2) \cong \S \gtimes C_0(E^{21}, \Cliff(p_{21}^*E^{21})) \gtimes_{_{E^2}} p_{21}^*C_0(E^1, \Cliff(p_1^*E^1)) \gtimes_{_M} \CC(M).
\end{eqnarray}
Note that since since $p_{21}^*C_0(E^1, \Cliff(p_1^*E^1)) = C_0(E^2) \gtimes_{E^1} C_0(E^1, \Cliff(p_1^*E^1))$
\newline this implies that
\begin{eqnarray}\label{eqn:c*decomp2}
\A(E^2) & \cong & \S \gtimes C_0(E^{21}, \Cliff(p_{21}^*E^{21})) \gtimes_{E^1} C_0(E^1, \Cliff(p_1^*E^1)) \gtimes_{_M} \CC(M) \\
\notag   & \cong & \S \gtimes C_0(E^2, \Cliff(p_{21}^*E^{21})) \gtimes_{E^1} \CC(E^1) \cong \A(E^{21}) 
\end{eqnarray}
which is the decomposition for the Thom map $\Psi_{21} : \A(E^1) \to \A(E^2)$.

Let $C_2 = C_{E^2} : E^2 \to \Cliff(p_2^*E^2)$ be the Thom operator for $p_2 : E_2 \to M$ and let $C_1 = C_{E^1} : E^1 \to \Cliff(p_1^*E^1)$ be the
Thom operator for $p_1 : E^1 \to M$ as in Theorem \ref{thm:cliffop}. Let
$C_{21} : E^{21}\to \Cliff(p_{21}^*E^{21})$ be the Thom operator for the bundle $p_{21} : E^{21} \to E^1$. If $e = e^\perp \oplus e^1 \in E^2 =
E^\perp \oplus E^1$, then one computes:
$$C_2(e) = C_{p_{2}(e)}(e)  = C_{p_{21}(e)}(e^\perp) \gtimes 1  + 1 \gtimes C_{p_1p_{21}(e)}(e^1) $$ with
respect to the decomposition (\ref{eqn:cliffbdl2}) above. Thus, it implies that we have a tensor product decomposition of essentially self-adjoint unbounded
multipliers
$$C_2 = C_{21} \gtimes 1 + 1 \gtimes p_{21}^*C_1$$ 
induced by the $\cs$-algebra decomposition above.

It suffices to compute the composition $\Psi_{21} \circ \Psi_1$ on generators of $\A(M) = \S \gtimes \CC(M)$ of the form $\exp(-x^2) \gtimes s$ and $ x
\exp(-x^2) \gtimes s$ where $s \in \CC(M)$. We obtain (using isomorphisms (\ref{eqn:c*decomp}) and (\ref{eqn:c*decomp2})) the elements,
respectively\footnote{We drop the subscripts on tensor products for clarity.},
$$\exp(-x^2) \gtimes \exp(-C_{21}^2) \gtimes p_{21}^*\exp(-C_1^2) \gtimes s$$ and
\begin{eqnarray*} x \exp(-x^2)  \gtimes \exp(-C_{21}^2) \gtimes p_{21}^* \exp(-C_1^2) \gtimes s \\
 + \exp(-x^2) \gtimes C_{21} \exp(-C_{21}^2) \gtimes p_{21}^* \exp(-C_1^2) \gtimes s \\
 + \exp(-x^2) \gtimes \exp(-C_{21}^2) \gtimes p_{21}^*C_1 \exp(-C_1^2) \gtimes s
\end{eqnarray*} 
However, by Corollary \ref{cor:exp}, we have equalities  $$\exp(-C_2^2) = \exp(-C_{21})^2 \gtimes p_{21}^*\exp(-C_1^2)$$ and also 
$$C_2 \exp(-C_2^2) = C_{21} \exp(-C_{21}^2) \gtimes p_{21}^*\exp(-C_1^2) + \exp(-C_{21}^2) \gtimes p_{21}^* (C_1 \exp(-C_1^2)).$$ But under the map $\Psi_2$
the generators $\exp(-x^2) \gtimes s$ and $ x \exp(-x^2) \gtimes s$ are mapped, respectively, to the elements
$$\exp(-x^2) \gtimes \exp(-C_2^2) \gtimes s$$ 
and
$$x \exp(-x^2) \gtimes \exp(-C_2^2) \gtimes s + \exp(-x^2) \gtimes C_2\exp(-C_2^2) \gtimes s.$$ Therefore, the composition $\Psi_{21} \circ
\Psi_1$ agrees with $\Psi_2$ on generators since $p_2 = p_1 \circ p_{21} $, and hence the two maps are equal.
\EPf

Let $\E$ be a complete Euclidean space of infinite dimension. Denote by $O(\E)$ the orthogonal group of all invertible bounded linear maps $T : \E \to
\E$ which preserve the inner product:
$$\langle T(v), T(w) \rangle = \langle v, w \rangle, \quad \forall v, w \in \E.$$
Let $\EE$ be an infinite dimensional manifold \cite{AMR88,Lan62} modelled on $\E$.

\begin{dfn} 
A smooth map $p : \EE \to M$ is a {\bf Euclidean bundle modelled on $\E$} if it defines a smooth locally trivial fiber bundle $\EE$ with
fiber $\E$ and structure group $O(\E)$ in the norm $($or strong operator$)$ topology \cite{Lan62,Phi89}. 
\end{dfn}

It follows that $\EE$ is equipped with a metric $\metric_x$ on each fiber $\EE_x$, which varies smoothly in $x \in M$. That is, the map $(e, f)
\mapsto \langle e, f \rangle_{p(e)}$ is smooth on the fiber product $\EE \times_{_M} \EE = \{ (e, f) | p(e) = p(f)\}$.

\begin{dfn}\label{def:finsubbdl}
A subset $S^a \subset \EE$ will be called a {\bf subbundle} if there is a real vector bundle $p_a : E^a \to M$ and an
exact sequence 
$$\xymatrix{ 0 \ar[r] & E^a \ar[r]^{i} & \EE, }$$
of vector bundles such that $i(E^a) = S^a$. That is, $i : E^a \to \EE$ is a vector bundle morphism such that for each $x \in M$ the continuous linear map
$i_x : E^a_x \to \EE_x$ on fibers is split injective (i.e., with split range). This gives
$S^a$ the $($unique$)$ structure of a  vector bundle, and we will identify $E^a \cong S^a \subseteq \EE$. If $\rank(E^a) = n < \infty$
then we will call $E^a \subset \EE$ a {\bf finite rank subbundle}.
\end{dfn}

The following is adapted from Proposition 3.4.18 \cite{AMR88}.

\begin{prop} 
Let $p : \EE \to M$ be a Euclidean bundle modelled on $\E$ and let $q : E \to M$ be a finite rank vector bundle. If $f : E \to \EE$ is a vector
bundle map let $f_x : E_x \to \EE_x$ be the linear restriction to the fibers over $x \in M$ and let $\range(f) =
\bigcup_{x \in M} \range(f_x)$. Then $\range(f) \subset \EE$ is a finite rank subbundle if and only if $x \mapsto \rank(f_x)$ is locally constant on $M$.
\end{prop}

\begin{thm} 
Suppose $\E$ has countably infinite dimension. If $p : \EE \to M$ is a Euclidean bundle modelled on $\E$, there is an
increasing sequence 
$$M = E^0 \subset E^1 \subset E^2 \subset \cdots \subset E^n \subset \cdots \subset \EE$$ 
of finite rank subbundles such that $\rank(E^n) = n$ and $\EE = \bigcup_0^\infty E^n.$
\end{thm}

\Pf The locally trivial fiber bundles with structure group $O(\E)$ are classified by the cocycles in the sheaf cohomology $H^1(M, O(\E))$ (which is not a
group since $O(\E)$ is not abelian) \cite{RW98}. However, by Kuiper's Theorem \cite{Kui65}, the group $O(\E)$ is contractible in the norm (or strong
operator) topology. Hence, there is only the trivial cocycle corresponding to $\FF = M \times \E \to M$, the trivial Euclidean bundle over $M$ modelled
on $\E$. Thus, there is a vector bundle isomorphism $\Phi : \FF \to \EE$. Since $\E$ has countably infinite dimension, there is an
increasing sequence 
$$\{0\}= V^0 \subset V^1 \subset V^2 \subset \cdots \subset \E$$ 
of subspaces of $\E$ such that $\dim(V^n) = n$. Let $F^n = M \times V^n \subset \FF$ be
the associated finite rank subbundle of $\FF$ for all $n \ge 0$. It follows that if we define $E^n = \Phi(F^n) \subset \EE$, then the sequence
$\{E^n\}_0^\infty$ is our desired increasing sequence of finite rank subbundles. \EPf

\noindent {\bf Important Note:} Although every Euclidean bundle $\EE$ with infinite dimensional fiber is trivializable, the finite dimensional subbundles are not necessarily
trivializable. Moreover, if a group $G$ is acting on $E$ and $M$ in a compatible way, as we will consider in the next
section, then $\EE$ itself may not be {\it equivariantly} trivializable.

\begin{dfn} 
Let $E^a$ and $E^b$ be finite rank subbundles of $\EE$. We will define $E^a \preceq E^b$ if $E^a \subseteq E^b$ as subsets of $\EE$ and there is
an exact sequence $0 \to E^a \to E^b$ of vector bundles  such that the following diagram commutes:
$$\xymatrix{
 & 0 \ar[d] &  \\
0 \ar[r] & E^a \ar[d] \ar[r] & \EE \ar[d]^-{\id} \\
0 \ar[r] & E^b \ar[r] & \EE 
}$$
That is, $E^a$ is a subbundle of $E^b$ in a way that is compatible with their subbundle structures from $\EE$.
\end{dfn}

\begin{dfn} 
Let $p : \EE \to M$ be a Euclidean bundle bundle modelled on $\E$. Denote by
$\FB(\EE)$ the collection of all smooth finite rank subbundles $p_a : E^a \to M$. Note that $\FB(\EE)$ is a directed system under inclusion $\preceq$
of subbundles.
\end{dfn}

Given a finite rank subbundle $E^a \subset \EE$ there is an induced Euclidean metric $\metric^a$, which on the fibers $E^a_x$ is  given by
restricting the metric $\metric_x$ of the fiber $\EE_x$ to the subspace $E^a_x \subset \EE_x$. We also need a compatible connection $\nabla^a$
on $E^a$, and this can be done by equipping $\EE \to M$ with a compatible connection.

\begin{lem} 
Let $p : \EE \to M$ be a Euclidean bundle modelled on $\E$.  There exists a connection $\nabla : C^\infty(M, \EE) \to C^\infty(M , T^*M \otimes
\EE)$ which is compatible with the Euclidean metric $\metric$ on $\EE$.
\end{lem}

\Pf Standard partition of unity argument. \EPf

If $E^a \subset \EE$ is a finite rank subbundle, then there is an orthogonal decomposition of Euclidean bundles $\EE = E^\perp \oplus E^a$. Let
$p^\perp_a : \EE \to E^a$ denote the orthogonal projection onto $E^a$. We can then induce (project) a connection $\nabla^a$ on $E^a$ by the diagram
$$\xymatrix{  C^\infty(M, \EE) \ar[r]^-{\nabla} & C^\infty(M, T^*M \otimes \EE) \ar[d]^-{(1 \otimes p^\perp_{a})_*} \\
 C^\infty(M, E^a) \ar[u]^-{i_*} \ar[r]^-{\nabla^a} & C^\infty(M, T^*M \otimes E^a)
}$$ 
One can then easily check that $\nabla^a$ is compatible with the induced metric $\metric^a$.

Given any two finite rank subbundles $E^a \preceq E^b$ of $\EE$, we then have a canonical orthogonal decomposition
$E^b \cong E^\perp_{ab} \oplus E^a$
Thus, we can define a bundle projection $p_{ab} : E^b \to E^a$. Let $E^{ba}$ denote the vector bundle $p_{ab} : E^b \to E^a$. By the
discussion prior to Proposition \ref{prop:transitive}, this has a Euclidean bundle structure given by $E^{ba} \cong p_a^*E^\perp_{ab}$ and compatible
connection $\nabla^{ab} = p_a^*\nabla^\perp$ where $\nabla^b = \nabla^\perp \oplus \nabla^a$. By Theorem \ref{thm:Thomhom} there is a Thom
$*$-homomorphism 
$$\Psi_{ab} : \A(E^a) \to \A(E^b)$$
which will be the connecting map in our construction of the direct limit $\cs$-algebra associated to $\EE$.

\begin{lem}\label{lem:functoriality} 
Let $p : \EE \to M$ be an affine Euclidean bundle modelled on $\E$. If $E^a \preceq E^b
\preceq E^c$ are finite rank affine subbundles of $\EE$, then there is a commutative diagram
of Thom $*$-homomorphisms:
$$\xymatrix{
 &   \A(E^c) \\
\A(E^b) \ar[ur]^{\Psi_{bc}} & \\
 &  \A(E^a) \ar[ul]^{\Psi_{ab}} \ar[uu]_{\Psi_{ac}}
}$$
\end{lem}

\Pf Note that since $E^a$, $E^b$ and $E^c$ are fibered as Euclidean bundles with compatible connections over the Riemannian manifold $M$, they have
canonical Riemannian structures, as {\it manifolds}, from Lemma \ref{lem:splitting}. We also have a commuting diagram of bundle projections
$$\xymatrix{
 & E^c \ar[dl]_{p_{bc}} \ar[dd]^{p_{ac}} \\
E^b \ar[dr]_{p_{ab}} & \\
& E^a
}$$
Note that the Riemannian structure that $E^c$ has as a bundle over $E^a$ and $E^b$ is compatible with the structure coming from $E^c \to M$. Similarly for
$E^b$. The result now follows from Proposition \ref{prop:transitive} 
\EPf

\begin{dfn} 
Let $p : \EE \to M$ be a smooth affine Euclidean bundle modelled on $\E$. We define $\A(\EE)$ to be the direct limit
$\cs$-algebra
$$\A(\EE) = \vlim{a} \A(E^a),$$ 
where the direct limit is taken over the directed system $\FB(\EE)$ of all finite rank subbundles $E^a \subset \EE$,
using the $*$-homomorphisms
$\Psi_{ab} : \A(E^a) \to \A(E^b)$ above.
\end{dfn}

Note that since $\EE$ is infinite dimensional, $\A(\EE)$ is {\it not} the suspension of a $\cs$-algebra. However, $\A(\EE)$ is a graded, separable,
nuclear $\cs$-algebra.

By transitivity again, we have that the $*$-homomorphisms $$\Psi_a : \A(M) \to \A(E^a)$$ are compatible with the direct limit, so there is a canonical
Thom map
$$\Psi_p : \A(M) \to \vlim{a} \A(E^a) = \A(\EE).$$
Note that this is equivalent to identifying $M$ as the zero finite rank subbundle (section) of $\EE \to M$.

\begin{lem}\label{lem:bottmap}
If $M = \{0\}$ is a point, so that $\EE = \E \to \{0\}$, then 
$$\Psi = \beta : \A(0) \to \A(\E)$$
is the Bott homomorphism of Higson-Kasparov-Trout \cite{HKT98}
\end{lem}

By the finite rank Thom isomorphism theorem \ref{thm:Thomfrc}, the connecting maps $\Psi_{ab} : \A(E^a) \to \A(E^b)$ induce isomorphisms in $K$-theory:
$$\Psi_{ab*} : K_*(\A(E^a)) \to K_*(\A(E^b)).$$
Using the continuity of operator $K$-theory with respect to direct limits \cite{WO93}, we obtain:
$$K_*(\A(\EE)) = K_*(\vlim{a} \A(E^a)) \cong \vlim{a} K_*(\A(E^a)) \cong K_*(\A(M)).$$
Therefore, we have proved the non-equivariant

\begin{ThomThm} 
Let $\EE$ be an affine Euclidean bundle on $M$ modelled on $\E$.  The inclusion $M \subset \EE$ as the zero
subbundle induces an isomorphism in $K$-theory:
$$\Psi_* : K_*(\A(M)) \to K_*(\A(\EE)).$$
\end{ThomThm}

\section{The Equivariant Thom Isomorphism}

In order to formulate the equivariant version of our theorem, we need to choose
an appropriate class of group actions. The most popular class of group actions
on manifolds are proper actions because of their application to $K$-theory \cite{Phi89}, the index theory of elliptic operators \cite{Kas83,
Pat00}, and Baum-Connes Conjectures \cite{BCH94}. Since our theorem must reduce to the Bott periodicity
theorem \cite{HKT98} when $M = \{\bullet\}$ is a point, we must allow for the 
trivial action of a non-compact group on a point. But, non-compact groups do not act {\it properly} on a point.

Let $G$ be a smooth, second countable, locally compact, Hausdorff topological group (e.g., a countable discrete group.) A locally compact
space $M$ is a {\it $G$-space} if it has a continuous (or smooth if $M$ is smooth) $G$-action
$$G \times M \to M.$$
If $H$ is a closed subgroup of $G$, and $x \in M$ then an {\it $H$-slice} through $x$ is
an $H$-invariant open neighborhood $U$ of $x$ such that the map $G \times U \to M$ given by $(g, u) \mapsto gu$ descends to a $G$-equivariant homeomorphism 
$$G \times_H U  = (G \times U)/H \to GU$$ 
and such that $GU$ is an open neighborhood of $x \in M$. If $M$ is a smooth manifold on which $G$ acts smoothly via
diffeomorphisms, then we require that this be a $G$-equivariant diffeomorphism.

\begin{dfn}\label{def:isoreductive}
Let $(M, d)$ be a locally compact metric $G$-space. The action of $G$ on $M$ is {\bf isoreductive} if:
\ilist
\item  $G$ acts by isometries of $(M, d)$;
\item  for every $x \in M$ there is a closed subgroup $H$ of $G$ and an $H$-slice $U$ through $x$ which is $H$-equivariantly contractible.
\endilist
\end{dfn}

Note that if $M = \{\bullet\}$ is a point, the trivial $G$-action is vacuously isoreductive since we can choose $H = G$ and $G \times_G \{\bullet\} =
\{\bullet\}$. If $G$ is a countable discrete group and $M$ is a discrete $G$-space, then $M$ is isoreductive with respect to the trivial
metric $d(x, y) = 1$ ($x \neq y$) since we can choose $U = \{x\} = B(x, 1)$, $H = G_x$ (the isotropy subgroup of $x$) and $G
\times_{G_x} \{x\} \cong G\cdot x$ is the orbit of $x$. 

Recall that an action of $G$ on $M$ is called {\it proper} \cite{Abel74,Phi89} if the structural map
$$\aligned
G \times M & \to M \times M \\
(g, m) & \mapsto (gm, m)
\endaligned$$
is a proper map, i.e., the inverse image of a compact set is compact.

For the following, see the discussion in Example (1.4) \cite{BCH94}.

\begin{lem} 
If $G$ is a countable discrete group, every smooth isometric proper action of $G$ on a Riemannian manifold is isoreductive.
\end{lem}

Indeed, using a $G$-partition of unity \cite{Phi89} we can average any metric on $M$ to obtain a $G$-invariant metric. Also, for any $x \in M$ we can choose
$H = G_x$ the (finite) isotropy subgroup of $x$ and $G_x$-slice consisting of an open ball of small enough radius.

Recall that a $\cs$-algebra $A$ is called a {\it $G$-$\cs$-algebra} if $G$ acts continuously on $A$ by $\cs$-algebra automorphisms (i.e., for every $a \in A$ the map $G \to
A$ given by $g \mapsto ga$ is continuous.)

\begin{lem}\label{lem:isom}
An isometric action of $G$ on $M$ induces a continuous action of $G$ on $\A(M)$
as $\Z_2$-graded $\cs$-algebra automorphisms.
\end{lem}

\Pf Given the isometry $\ga : M \to M$ there is an induced Euclidean bundle isomorphism $T\ga : TM \to TM$
of the tangent bundle $TM$. This induces an isomorphism $\widehat{T\ga} : \Cliff(TM) \to \Cliff(TM)$
of the $\Z_2$-graded Clifford bundle. Hence, we obtain an induced automorphism $\ga_\tau : \CC(M) \to \CC(M)$ of $\Z_2$-graded $\cs$-algebras
given by $(\ga_\tau s)(x) = \widehat{T\ga}(s(\ga^{-1}x))$ for all $x \in M$ and $s \in \CC(M)$. Tensoring with $\id_\S$
gives the automorphism $\ga_* = \id_\S \gtimes \ga_\tau : \S \CC(M) =  \A(M) \to \A(M)$. \EPf

\begin{dfn}\label{def:ktheory} 
If $A$ is a $G$-$\cs$-algebra, we define the {\it equivariant $K$-theory} $K^G_*(A)$ to
the be the $K$-theory of the full crossed product $\cs$-algebra $A \rtimes
G$: $$K^G_j(A) = K_j(A \rtimes G), \quad j = 0, 1.$$
\end{dfn}

For compact groups, the above definition is actually the Green-Julg-Rosenberg Theorem 11.7.1 \cite{Bla98}.
See Pedersen \cite{Ped79} for a discussion of full crossed product
$\cs$-algebras. It follows that if $\phi : A \to B$ is an equivariant
$*$-homomorphism then there is an induced $\cs$-algebra $*$-homomorphism
$\phi^G : A \rtimes G \to B \rtimes G$. Hence, there is an induced
homomorphism  $\phi_*^G : K^G_*(A) \to K^G_*(B)$ of equivariant $K$-theory groups.
We chose the full crossed product because of this universal property and
its compatibility with equivariant $E$-theory \cite{GHT00}.

Note that if the action of $G$ on $M$ is proper then $K^*_G(M)$, the equivariant topological $K$-theory of $M$ (as defined by N.C. Phillips \cite{Phi89}) exists. This
$K$-theory satisfies a generalized Green-Julg-Rosenberg theorem:
$$K^*_G(M) \cong K_*(C_0(M) \rtimes_\alpha G)$$
where $C_0(M) \rtimes_\alpha G$ is {\it any} (e.g., full or reduced) crossed product of $C_0(M)$ by $G$. 

\begin{prop} 
If $M$ is an oriented even-dimensional \spinc-manifold and $G$ acts on $M$ properly such that the \spinc-structure is $G$-invariant, then
there are canonical isomorphisms
$$K^G_*(\A(M)) \cong K^{*+1}_G(M)$$
in equivariant $K$-theory.
\end{prop}

From now on, we assume that $M$ is a smooth finite-dimensional $G$-manifold equipped with a fixed $G$-invariant Riemannian metric and the
action of $G$ is isoreductive with respect to this metric.

\begin{dfn} 
Let $p_a : E^a \to M$ and $p_b : E^b \to M$ be smooth Euclidean bundles on $M$. 
If $\ga : M \to M$ is an isometry, then a bundle map $\gg : E^a \to E^b$ is called a 
{\bf Euclidean bundle isomorphism over} $\ga$  if the following diagram commutes
$$\xymatrix{
E^a \ar[r]^-{\gg} \ar[d]_{p_a} & E^b \ar[d]^{p_b} \\
M \ar[r]^-{\ga} & M 
}$$
and $\gg_x : E^a_x \to E^b_{\ga(x)}$ is an isomorphism of Euclidean spaces for each $x \in M$.
\end{dfn}

Let $\E$ be a Euclidean space of countably infinite dimension.

\begin{dfn} 
A {\bf Euclidean $G$-bundle} $\EE$ on $M$ with fiber $\E$ is a smooth
Euclidean bundle $p : \EE \to M$ modelled on $\E$ equipped with a smooth 
action $G \times \EE \to \EE$ by Euclidean bundle automorphisms such that
the projection $p : \EE \to M$ is $G$-equivariant. That is, for each $\ga \in G$ there is
a Euclidean bundle automorphism $\gg : \EE \to \EE$ over $\ga : M \to M$
such that $\gg_1 \circ \gg_2 = \gg$ where $\ga = \ga_1 \ga_2 \in G$.
\end{dfn}

\begin{lem}\label{lem:g-subbdl}
Let $p : \EE \to M$ be a Euclidean $G$-bundle with fiber $\E$. Then for any
smooth finite rank subbundle $E^a \subset \EE$ and $\ga \in G$, the
image $\ga \cdot E^a$ has a canonical structure as a finite rank subbundle and 
$$\gg_! = \gg|_{E^a} : E^a \to \ga \cdot E^a$$
is a finite rank Euclidean bundle isomorphism over $\ga : M \to M$ .
\end{lem}

\Pf By Definition \ref{def:finsubbdl}, there is a subset $S^a \subset \EE$, a finite rank Euclidean bundle
$p_a : E^a \to M$ and an exact sequence of Euclidean bundles $$\xymatrix{ 0 \ar[r] & E^a \ar[r]^{i_a} & \EE }$$
such that $i(E^a) = S^a$. Now, define the subbundle $\ga \cdot E^a \subset \EE$ as follows. As a subset of $\EE$ we have
$\ga \cdot E^a = \ga \cdot i_a(E^a) = \ga \cdot S^a = \{\ga \cdot e : e \in S^a\}$. As a finite rank Euclidean bundle we
use $p_a : E^a \to M$ but with associated exact sequence 
$$\xymatrix{ 0 \ar[r] & E^a \ar[r]^{\gg \circ i_a} & \EE }$$ 
given by the composition
$$\xymatrix{ 0 \ar[r] & E^a \ar[r]^{i_a} & \EE \ar[r]^-{\gg} & \EE }$$
It follows that there is a commutative diagram 
$$\xymatrix{
0 \ar[r] & E^a \ar[r]^-{i_a} \ar[d]_{\gg} & \EE \ar[d]^-{\gg} \\
0 \ar[r] &  E^a \ar[r]^-{\gg \circ i_a} & \EE
}$$
where $\gg_! : E^a \to E^a$ is the Euclidean bundle automorphism $i_a^{-1} \circ \gg  \circ i_a : E^a \to E^a$.
\EPf

\begin{rmk}
We should make the important observation that the action of $G$ on $\EE$ may not preserve {\bf any} finite rank subbundle except, of course, the zero
subbundle $M \subset \EE$. Also, if the action of $G$ on $M$ has no fixed points, then the fibers $\EE_x \cong \E$ may carry no action of $G$.
\end{rmk}

\begin{dfn} 
Let $(E^a, \nabla^a)$ and $(E^b, \nabla^b)$ be Euclidean bundles on $M$ with compatible connections.
Let $\ga : M \to M$ be an isometry of $M$. If $\gg : E^a \to E^b$ is a Euclidean bundle isomorphism
over $\ga$, then we call $\gg$ an {\bf affine Euclidean isomorphism} $($over $\ga )$ if the following diagram commutes:
$$\xymatrix{
C^\infty(M, E^b) \ar[r]^-{\nabla^b} & C^\infty(M, T^*M \otimes E^b) \\
C^\infty(M, E^a) \ar[r]^-{\nabla^a} \ar[u]^{\gg_*} & C^\infty(M, T^*M \otimes E^b) \ar[u]_{(T^*\ga \otimes \gg)_*}
}$$
where $T^*\ga : T^*M \to T^*M$ is the induced Euclidean bundle isomorphism of the cotangent bundle $T^*M \cong TM$.
We denote this by $(\gg, \ga) : (E^a, \nabla^a) \to (E^b, \nabla^b)$.
\end{dfn}

\begin{lem} 
Let $(\gg, \ga) : (E^a, \nabla^a) \to (E^b, \nabla^b)$ be an affine Euclidean isomorphism
of finite rank affine Euclidean bundles on $M$. There is an induced isomorphism
of $\Z_2$-graded $\cs$-algebra $(\gg, \ga)_* : \A(E^a) \to \A(E^b)$ such that the following
diagram commutes
$$\xymatrix{
\A(E^a) \ar[r]^-{(\gg, \ga)_*}_{\cong} &  \A(E^b) \\
\A(M) \ar[u]^-{\Psi_a} \ar[r]^{\ga_*}_{\cong} & \A(M) \ar[u]_-{\Psi_b}
}$$
where $\ga_*$ is the isomorphism from Lemma \ref{lem:isom} and $\Psi_a$ and $\Psi_b$ are the Thom maps from Theorem \ref{thm:Thomhom}.
\end{lem}

\Pf 
Using the connections and the fact that $(\gg, \ga)$ transforms the connection $\nabla^a$ into the connection $\nabla^b$,
there is a commutative diagram
$$\xymatrix{
TE^b  =              & p_b^*E^b \oplus p_b^*TM \\
TE^a = \ar[u]^-{T\gg}  & p_b^*E^a \oplus p_a^*TM \ar[u]^-{\gg^* \oplus T\ga^*}_-{\cong}
}$$
where the vertical maps appearing on the RHS are the Euclidean bundle isomorphisms induced by the pair $(\gg, \ga)$. Thus, with the metrics
induced by the splitting, the induced map $T\gg : TE^a \to TE^b$ is a Euclidean bundle isomorphism, which implies that
$\gg : E^a \to E^b$ induces an {\it isometry} between the Riemannian manifolds $E^a$ and $E^b$.

There is also an induced commutative diagram
$$\xymatrix{
E^b \ar[r]^-{\tau_b} \ar@/^2pc/[rr]^-{C_b}                   & p_b^*E^b \ar[r]^-{i_b}                        & \Cliff(p_b^*E^b) \\
E^a \ar[u]^-{\gg}_-{\cong} \ar[r]^-{i_b} \ar@/_2pc/[rr]^-{C_a} & p_a^*E^a \ar[r]^-{i_b} \ar[u]^-{\gg^*}_-{\cong} & \Cliff(p_a^*E^a) \ar[u]^-{\gg^*}_-{\cong}
}$$
where the vertical maps are the appropriate bundle isomorphisms induced by the map $\gg$.
Thus, the isomorphism $\gg$ transforms the Thom operator $C_a$ of $E^a$ into the Thom operator $C_b$ of $E^b$ in a manner
compatible with the isometric action of $\ga : M \to M$. The result now easily follows using Theorems \ref{thm:decomp}
and \ref{thm:Thomhom} and the equivariant properties of the functional calculus.
\EPf

\begin{dfn}
Let $p : \EE \to M$ be a Euclidean $G$-bundle. A connection $\nabla$ on $\EE$ is called {\bf $G$-invariant}
if for each $\ga \in G$ the bundle map $\gg : \EE \to \EE$ determines an affine Euclidean isomorphism
$(\gg, \ga) : (\EE, \nabla) \to (\EE, \nabla).$
That is, $$(T^*\ga \otimes \gg)_* \circ \nabla =  \nabla \circ \gg_*$$ 
for all $\ga \in G$. If $\EE$ is equipped with a $G$-invariant connection, then we say that $\EE$ is
an {\bf affine} Euclidean $G$-bundle.
\end{dfn}

Averaging any connection on $\EE$ with respect to a $G$-partition of unity yields

\begin{lem} 
If the action of $G$ on $M$ is proper, there exists a $G$-invariant connection $\nabla$ on any Euclidean $G$-bundle $p : \EE \to M$.
\end{lem}

\begin{lem} 
Let $p : \EE \to M$ be an affine Euclidean $G$-bundle with $G$-invariant connection $\nabla$. Let $E^a \subset \EE$ be a finite rank subbundle.
The map $\gg_! : E^a \to \ga \cdot E^a$ $($from Lemma \ref{lem:g-subbdl}$)$ determines an affine Euclidean bundle isomorphism
$$(\gg_!, \ga) : (E^a, \nabla^a) \to (\ga \cdot E^a, \ga \cdot \nabla^a)$$
where $\nabla^a$ denotes the connection induced on $E^a$ by $\nabla$.
\end{lem}

\begin{cor}\label{cor:g-act}
With the above hypotheses on $(\EE, \nabla, M, G)$, we have that for each finite rank subbundle $E^a \subset \EE$ and each $\ga \in G$
there is a commutative diagram
$$\xymatrix{
\A(E^b) \ar[r]^-{(\gg_!, \ga)_*}_-{\cong} &  \A(\ga \cdot E^a) \\
\A(M) \ar[u]^-{\Psi_a} \ar[r]^-{\ga_*}_-{\cong} & \A(M) \ar[u]_-{\ga \cdot \Psi_a}
}$$
where $\ga \cdot \Psi_a : \A(M) \to \A(\ga \cdot E^a)$ is the Thom map for the finite rank subbundle $\ga \cdot E^a \subset \EE$.
\end{cor}

\begin{lem}\label{lem:bundleisom}
With the above hypotheses on $(\EE, \nabla, M, G)$, there is an induced continuous action of $G$ on the direct limit $\cs$-algebra
$$\A(\EE) = \vlim{a} \A(E^a)$$
as $\Z_2$-graded $\cs$-algebra automorphisms.
\end{lem}

Moreover, the $G$-equivariant inclusion $M \subset \EE$ as the zero finite rank subbundle induces a $G$-equivariant
$\Z_2$-graded $*$-homomorphism
$$\Psi_p : \A(M) \to \A(\EE).$$
As discussed in the last section, we know that $\Psi_p$ induces an isomorphism on $K$-theory, but we want to show
furthermore that $\Psi_p$ induces an isomorphism in {\it equivariant} $K$-theory. 
The main result of this paper is as follows.

\begin{thm}\label{GThomThm}
Let $G$ be a smooth, second countable, locally compact group and $M$ be a smooth isoreductive Riemannian $G$-manifold. If $p : \EE \to M$ is an affine
Euclidean $G$-bundle with fiber a countably infinite dimensional Euclidean space 
then the inclusion of $M \subset \EE$ as the zero finite rank subbundle induces an isomorphism in equivariant $K$-theory:
$$\Psi_*^G : K^G_*(\A(M)) \to K^G_*(\A(\EE)).$$
\end{thm}

If the Euclidean space $\E$ is equipped with an orthogonal action of $G$, then we may
consider the trivial Euclidean $G$-bundle $p: \E \to \{0\}$, where $M = \{0\}$ denotes the isoreductive (but very
non-proper!) one-point Riemannian $G$-manifold consisting of the origin. Our theorem then reduces
to the Bott periodicity theorem. The Thom map $\Psi_p = \beta$ is the ``Bott map'' in this case.

\section{Proof of the Thom Isomorphism Theorem \ref{GThomThm}}

Let $G$ be a smooth, second countable, locally compact, Hausdorff group.

\begin{dfn}\label{def:isoreductalg} 
A $G$-$\cs$-algebra $A$ is called {\bf isoreductive} if there is a a second countable, locally compact,
isoreductive metric $G$-space $(M, d)$ and a nondegenerate $G$-equivariant injection 
$$\theta : C_0(M) \to ZM(A)$$ 
such that $\theta(C_0(M))A$ is dense in $A$. If $f \in C_0(M)$ and $a \in A$ we write $fa = \theta(f)a$.
If $A$ is $\Z_2$-graded then we assume that $\theta$ has grading degree zero, where  $C_0(M)$ is
trivially graded. That is, $A$ is a $\Z_2$-graded $C_0(M)$-algebra $($Definition \ref{def:grddCo(M)alg}$)$ such that the structural
homomorphism is $G$-equivariant.
\end{dfn}

If $A$ is an isoreductive $G$-$\cs$-algebra, and $M$ is the isoreductive metric space in the previous definition, then we will say that
$A$ is {\it isoreductive over} $M$.  For a discussion of the related notion of (ungraded) proper algebras, see \cite{GHT00}.

\begin{lem} 
Let $A$ be a $\Z_2$-graded $G$-$\cs$-algebra. If $A$ is isoreductive over $M$ then $B \gtimes A$ is isoreductive over
$M$ for any $\Z_2$-graded $G$-$\cs$-algebra $B$.
\end{lem}

\begin{ex} 
Let $M$ be a second countable, isoreductive Riemannian $G$-manifold. The $\Z_2$-graded $G$-$\cs$-algebra $\A(M) = \S
\gtimes \CC(M)$ is isoreductive over $M$, where $\S$ has the trivial $G$-action.
\end{ex}

Let $\E$ be a  Euclidean space of countable dimension.

\begin{prop}\label{prop:isoreductiveSCE} 
Let $p : \EE \to M$ be an affine Euclidean $G$-bundle modelled on $\E$.
The $\cs$-algebra $\A(\EE)$ is isoreductive over $M$.
\end{prop}

\Pf If $E^a \subset \EE$ is a finite rank subbundle, then using the projection $p_a : E^a \to M$, $\A(E^a)$ has a graded
$C_0(M)$-algebra structure. Moreover, if $E^a \preceq E^b$ are finite rank subbundles, the Thom map $\A(E^a) \to \A(E^b)$ from
Theorem \ref{thm:Thomhom} is clearly $C_0(M)$-linear. It follows by an approximation argument that the direct limit
$$\A(\EE) = \vlim{a} \A(E^a)$$ has a graded $C_0(M)$-algebra structure and the induced action of $C_0(M)$ is $G$-equivariant. \EPf

If a $\cs$-algebra $A$ is isoreductive over $M$, then $A$ is a $C_0(M)$-algebra and can thus be realized as the algebra
of sections of an (upper-semicontinuous) $\cs$-bundle over $M$ as follows \cite{EW97,Nil96}. Given $x \in M$, the fiber over $x$ is the
quotient $A_x = A/I_x$ where $I_x$ is the ideal $I_x = C_0(M\minus\{x\})A$. For any $a \in A$, let $a(x)$ denote the image of $a$ 
in the fibre $A_x$. There is then a faithful representation of $A$ into the direct sum $\bigoplus_{x \in M}A_x$ given by $a \mapsto (a(x))_{x \in M}$.

\begin{dfn} $($Compare Definition 2.1\cite{EW97}$)$ 
Let $A$ be isoreductive over $M$. Let $Y$ be a locally compact subset of $M$. The {\bf restriction} of $A$ to $Y$ is defined as
$$A_Y = C_0(Y) \cdot A = \{b \in \oplus_{y\in Y}A_y: b(y) = f(y)a(y) \text{ for some } f \in C_0(Y) \text{ and } a \in A\}.$$
\end{dfn}
\noindent This definition makes sense for any $C_0(M)$-algebra (without grading or $G$-action.) The following is adapted
from Lemma 2.2 \cite{EW97}.

\begin{lem}\label{lem:restrictU}
Let $A$ be isoreductive over $M$. Let $Y$ be a $G$-invariant locally compact subset of $M$. Then $A_Y$ is a $G$-$C_0(Y)$-algebra.
If $U$ is a $G$-invariant open subset of $M$, then $A_U$ is a $G$-invariant ideal in $A = A_M$ and can be identified with $C_0(U)A$. Moreover,
there is an exact sequence
$$\xymatrix{ 0 \ar[r] & A_U \ar[r] & A \ar[r] & A_{M \minus U} \ar[r] & 0 }$$ 
of $G$-$\cs$-algebras.
\end{lem}

\begin{lem} 
Let $A$ be isoreductive over $M$. If $U$ is a $G$-invariant open subset, then there is a cyclic six term exact sequence
$$\xymatrix{
K^G_0(A_U) \ar[r] & K^G_0(A) \ar[r] & K^G_0(A_{M \minus U}) \ar[d]  \\ 
K^G_1(A_{M \minus U})  \ar[u]  & K^G_1(A) \ar[l] & K^G_1(A_U) \ar[l]
}$$
\end{lem}

\Pf  From the previous lemma, we have that
$$\xymatrix{ 0 \ar[r] & A_U \ar[r] & A \ar[r] & A_{M \minus U} \ar[r] & 0
}$$ is an exact sequence of $G$-$\cs$-algebras. Thus, by the functorial properties of the full crossed product, we have an
induced short exact sequence
$$\xymatrix{ 0 \ar[r] & A_U \rtimes G \ar[r] & A \rtimes G\ar[r] & A_{M \minus U} \rtimes G \ar[r] & 0
}$$ of crossed products. Now apply the $K$-theory six-term exact sequence~\cite{WO93}.
\EPf

The proof of the following result does not seem to be in the literature, but was communicated to the author by D. Williams.

\begin{lem} 
Let $A$ be a $G$-$\cs$-algebra. Suppose that $I$ and $J$ are $G$-invariant ideals.  Then $I \rtimes G$ and $J \rtimes G$
are ideals in $A \rtimes G$. Let $I \cap J$ denote the intersection ideal. Then,
\begin{equation}\label{eq:1}
  (I\rtimes G) \cap (J\rtimes G)=(I \cap J)\rtimes G.
\end{equation}
\end{lem}

\Pf
First recall that $I\rtimes G=\overline{C_{c}(G,I)}$ where the latter
is viewed as an ideal of $C_{c}(G,A)\subset A\rtimes G$.  Since 
\begin{equation*}
  C_{c}(G,I)\cap C_{c}(G,J)\subset C_{c}(G,I\cap J),
\end{equation*}
we certainly have
\begin{equation*}
  (I\rtimes G)\cap (J\rtimes G)\subset (I \cap J)\rtimes G.
\end{equation*}
On the other hand, $(I \cap J)\rtimes G$ is the closure of 
\begin{equation*}
  \operatorname{span}\{\,f\otimes c:\text{$f\in C_{c}(G)$ and $c\in I \cap J$}\,\},
\end{equation*}
where $f\otimes c$ is the shorthand for the function $g \mapsto f(g)c$.  
By the Cohen factorization theorem \cite{Co59} or an approximation argument, we can
assume that $c=ab$ with $a\in I$ and $b\in J$.  Now if $\{\,e_{i}\,\}$
is an approximate unit in $C_{c}(G)$, then 
\begin{equation*}
  (e_{i}\otimes a)(f\otimes b)=e_{i}*f\otimes c
\end{equation*}
converges to $f\otimes c$ in $(I\cap J)\rtimes G$.  It follows that we get
equality in \eqref{eq:1}. \EPf

The following is an equivariant Mayer-Vietoris theorem.

\begin{lem}\label{lem:G-mayer} 
Let $A$ be isoreductive over $M$. If $U$ and $V$ are $G$-invariant open subsets of $M$, there is an exact sequence
$$\xymatrix{
K_0^G(A_{U \cap V}) \ar[r] & K_0^G(A_U) \oplus K_0^G(A_V) \ar[r] & K_0^G(A_{U \cup V}) \ar[d] \\
K_1^G(A_{U \cup V}) \ar[u] & K_1^G(A_U) \oplus K_1^G(A_V) \ar[l] & K_1^G(A_{U \cap V}) \ar[l]
}$$
of equivariant $K$-theory groups. 
\end{lem}

\Pf Since $U$ and $G$ are $G$-invariant open sets, it follows that $A_U$, $A_V$, and $A_{U \cap V}$ are $G$-invariant ideals in $A_{U \cup V}$.
We also have the $G$-equivariant identifications
$$ \Bigg\{ \aligned 
A_U + A_V &= A_{U \cup V}\\
A_U \cap A_V &= A_{U\cap V}.
\endaligned$$
It follows, using the previous lemma, that
$$ \Bigg\{  \aligned
A_U \rtimes G + A_V \rtimes G & = (A_U +A_V) \rtimes G && = A_{U \cup V} \rtimes G\\
(A_U \rtimes G) \cap (A_V \rtimes G) & = (A_U \cap A_V) \rtimes G && = A_{U\cap V}\rtimes G.
\endaligned$$
The result now follows from Definition \ref{def:ktheory} and Exercise 4.10.21 \cite{HR00}. \EPf

\begin{lem}\label{lem:directlimU} 
Let $A$ be isoreductive over $M$.
Suppose $M = \bigcup_{i \in I} U_i$ where $\{U_i\}_{i \in I}$ is an increasing directed system of
$G$-invariant open subsets of $M$. That is, $U_i
\subseteq U_j$ if $i < j$ . Then we have isomorphisms
$$\aligned
A & \cong \vlim{i \in I} \A_{U_i} \\ 
K^G_*(A) & \cong \vlim{i \in I} K^G_*(A_{U_i}).
\endaligned$$
\end{lem}

\Pf The first direct limit isomorphism follows by an easy approximation argument. The $K$-theory isomorphism follows from the first and
the fact that the $A_{U_i}$ are $G$-invariant ideals, which implies that there is a direct limit isomorphism
$$A \rtimes G \cong \vlim{i} A_{U_i} \rtimes G$$ of (full) crossed product $\cs$-algebras. Now use the fact that $K$-theory is
continuous, i.e., commutes with direct limits.
\EPf

If $U \subset M$ is open, then for any finite rank subbundle $E^a \subset \EE$ we have that $E^a|_U \subset E^a$ as an open subset
and $E^a|_U \subset \EE|_U$ as a finite rank subbundle. Moreover, if $E^a \preceq E^b \preceq E^c$ then we have that
$E^a|_U \preceq E^b|_U \preceq E^c|_U$. Thus, we are led to make the following definition.

\begin{dfn} 
Let $p : \EE\to M$ be an affine Euclidean $G$-bundle. For any $(G$-invariant$)$ open subset
$ i : U \hookrightarrow M$ we define the $\cs$-algebra
$$\A(\EE|_U) =_{\text{def}}  \vlim{a} \A(E^a|_U)$$ where the direct limit is taken over all finite rank subbundles $E^a \subset \EE$.
Note that we give the restricted bundle
$\EE|_U = i^*(\EE)$ the connection $\nabla_U = i^*(\nabla)$.
\end{dfn}

\begin{lem}  
Let $M$ be an isoreductive Riemannian $G$-manifold. Let $\EE \to M$ be an affine Euclidean $G$-bundle modelled on $\E$. 
If $U$ is an open subset of $M$, there are isomorphisms
$$\aligned
& \A(M)_U && \cong \A(U) \\
& \A(\EE)_U && \cong \A(\EE|_U). 
\endaligned$$ If $U$ is $G$-invariant, these isomorphisms are $G$-equivariant.
\end{lem}

\Pf The first isomorphism follows from the definitions and the fact that $\CC(M)_U \linebreak[0] = C_0(U)\CC(M) = \CC(U)$. Given any finite rank subbundle $E^a
\subset \EE$, it is easy to see that all elements in
$\A(E^a|_U)$ are in $\A(\EE)_U$. This implies that the direct limit $\A(\EE|_U) \subseteq \A(\EE)_U$. For the reverse, it follows that
$\A(\EE)_U \cong \vlim{} \A(E^a)_U$, but
$\A(E^a)_U \subset \A(E^a|_U)$ and so $\A(\EE)_U \cong \A(\EE|_U)$ as desired. \EPf

Before we can begin the proof, we need the following definition \cite{RW98}. Let $H$ be a closed subgroup of $G$ and let $B$ be an
$H$-$\cs$-algebra. Note that if 
$f : G \to B$ satisfies $f(\ga h) = h^{-1}(f(\ga))$ for all $\ga \in G$ and $h \in H$, then the function $\ga \mapsto \|f(\ga)\|$ is
constant on cosets $\ga \cdot H$ and thus determines a well-defined function $\|f(\cdot H)\|$ on the quotient space $G/H$
\begin{dfn} Let $H$ be a closed subgroup of $G$ and let $B$ be an $H$-$\cs$-algebra. The {\bf induced $\cs$-algebra} $\Ind^G_H B$ is
defined as formula
$$\Ind^G_H B = \{f \in C_b(G, B) : f(gh) = h^{-1}(f(g)) \text{ \& } \|f(\cdot H)\| \in C_0(G/H)\}.$$ The group $G$ acts on $\Ind^G_H
B$ by left translation: $(\ga f)(g) = f(\ga^{-1}g)$. A $\Z_2$-grading on $B$ induces a $\Z_2$-grading on $\Ind^G_H B$ in the obvious
way.
\end{dfn}

For example, if $U$ is an $H$-space then $\Ind_H^G C_0(U) \cong C_0(G \times_H U)$.

\begin{lem}\label{lem:induct} 
Let $H$ be a closed subgroup of $G$. If $B$ is an $H$-$\cs$-algebra, there are canonical isomorphisms
$$K_*^H(B) \cong K_*^G(\Ind^G_H B).$$
\end{lem}

\Pf By Green's Symmetric Imprimitivity Theorem \cite{RR88}, $\Ind_H^G B \rtimes G$ and $B \rtimes H$ are strongly Morita equivalent.
Hence, they have the same $K$-theory and so:
$$K_*^H(B) = K_*(B \rtimes H) \cong K_*(\Ind_H^G B \rtimes G) = K_*^G(\Ind_H^G B).$$ \EPf

\begin{prop}\label{prop:induct}  
Let $A$ be isoreductive over $M$. For any $x \in M$, let $U$ be an $H$-slice through $x$ as in Definition (\ref{def:isoreductive}). 
There is a $G$-equivariant isomorphism
$$\Ind^G_H A_U  \cong  A_{G U}.$$
Hence, $K_*^H(A_U) \cong K_*^G(A_{GU})$.
\end{prop}

\Pf 
Since $U$ is $H$-invariant, it follows that the ideal $A_U$ is an $H$-$\cs$-algebra and has a $C_0(U)$-algebra structure. Similarly, $A_{GU}$ is
a $G$-$C_0(GU)$-algebra and isoreductive over $GU$. Thus, by Proposition 2.1 \cite{Nil96}, there is a canonical continuous $G$-equivariant map (with dense range)
$$\Res : \Prim A_{GU} \to GU$$
where $$\Res(J) = x \iff \overline{A_{GU}I_xA_{GU}} = C_0(GU\minus\{x\})A = A_{GU\minus\{x\}}\subset J$$ 
where $I_x = \{f \in C_0(M) : f(x) = 0\}$.
Since $U$ is an $H$-slice, there is a $G$-equivariant homeomorphism
$$\phi : GU \to G \times_H U.$$
Let $p : G \times_H U \to G/H$ be the $G$-map induced by the projection $G \times U \to G$.
Thus, we have a canonical continuous $G$-equivariant map
$$\sigma = p \circ \phi \circ \Res : \Prim A_{GU} \to G/H.$$
Let $I$ be the ideal 
$$I = \bigcap \{J \in \Prim A_{GU} : \sigma(J) = eH\}.$$
Note that $(p \circ \phi)^{-1}(eH) = U.$
One can then check that
$I = A_{GU \minus U}$.
By a theorem of Echterhoff \cite{Ech90,Ech92}, there is a
$G$-equivariant isomorphism
$$\Phi : A_{GU} \to \Ind_H^G A_{GU}/I.$$
Since $(p \circ \phi)^{-1}(eH) = U$, $U$ is a closed subset of $GU$. By Lemma \ref{lem:restrictU}: 
$$A_{GU}/I = A_{GU}/A_{GU \minus U} \cong A_U.$$ 
Thus, there is an induced isomorphism of crossed products
$$\Phi^G : A_{GU} \rtimes G \to \Ind_H^G A_U \rtimes G$$
and the result follows.
\EPf

Let $\E$ be a Euclidean space of countably infinite dimension (with orthogonal $G$-action.)

\begin{prop} \label{prop:trivisom}
Let $\EE = M \times \E \to M$ be a trivial affine Euclidean $G$-bundle modelled on $\E$ $($with the trivial $G$-connection$)$.
There is a $G$-equivariant isomorphism 
$$\A(\EE) \cong \A(\E) \gtimes \CC(M)$$ of $\Z_2$-graded $G$-$\cs$-algebras. Under this isomorphism, the Thom map has the form
$$\Psi = \beta \gtimes \id_{\CC(M)} : \S \gtimes \CC(M) \to \A(\E) \gtimes \CC(M)$$
where $\beta : \S = \A(0) \to \A(\E)$ is the Bott map of Lemma \ref{lem:bottmap}. Thus,
the induced map $\Psi^G_* : K^G_*(\A(M)) \to
K^G_*(\A(\EE))$ is an isomorphism on equivariant $K$-theory.
\end{prop}

\Pf Let $\{V_\alpha\}$ be the collection of all finite dimensional subspaces of $\E$ ordered by subspace inclusion. Then $F^\alpha = M
\times V_\alpha$ defines a directed system
$\{F_\alpha\}$ of finite rank Euclidean subbundles such that $\EE = \bigcup_\alpha F^\alpha$. By Proposition
\ref{prop:triv} we have that $$\A(F^\alpha) \cong \A(V_\alpha) \gtimes \CC(M).$$ Moreover, if $\alpha < \beta$ then the Thom map
has the form
$$\Psi_{\alpha \beta} = \beta_{\alpha \beta} \gtimes \id_{\CC(M)} : \A(V_\alpha) \gtimes \CC(M) \to \A(V_\beta) \gtimes \CC(M).$$
This implies that we have a
$*$-homomorphism
$$\theta : \A(\E) \gtimes \CC(M) = \vlim{\alpha} \A(F^\alpha) \to \A(\EE).$$ Given any finite rank subbundle $E^a \subset \EE$, there
is a $\beta$ such that $E^a \preceq F^\beta$. This implies that the  directed subsystem $\{F^\alpha\}$ is cofinal with respect to the
directed system $\FB(\EE) = \{E^a\}$. By Lemma L.1.5 \cite{WO93}, there is a
$*$-homomorphism 
$$\A(\EE) \to \vlim{\alpha} \A(F^\alpha) \cong \A(\E) \gtimes \CC(M)$$ which provides the inverse to $\theta$. If $E^a \preceq
F^\beta$ then for any $\ga \in G$, we have that
$\ga \cdot E^a \preceq \ga \cdot F^\beta = M \times \ga \cdot V_\beta = M \times V^\delta$ and so these isomorphisms are
$G$-equivariant. Using the equivariant asymptotic morphism calculations in Sections 4, 5 and Appendix C of \cite{HKT98} we easily obtain (by tensoring with
$\id_{\CC(M)}$) the last statement. \EPf

\Pff {\it of Theorem \ref{GThomThm}}. Using Definition \ref{def:isoreductive}, Mayer-Vietoris and Five Lemma arguments, direct limits, and transfinite induction,
it suffices to show that the Thom map induces an isomorphism for $W = GU$, where $U$ is an $H$-slice for some closed subgroup $H$ of $G$. 
Since $U$ is $H$-equivariantly contractible, we may also assume that the restricted bundle
$\EE|_U \cong U \times \E$ is a trivial affine Euclidean $H$-bundle (with the trivial
$H$-invariant connection.) By Proposition \ref{prop:trivisom},
$$\Psi^H_* : K^H_*(\A(U)) \to K^H_*(\A(\EE|_U))$$
is an isomorphism of abelian groups. Now, by Lemma \ref{lem:induct} and Proposition \ref{prop:induct}, it follows for $W=GU$ that
$$\aligned K_*^G(\A(W)) & \cong K_*^G(\Ind^G_H \A(U)) && \phantom{} \\
                   & \cong K_*^H(\A(U))      && \cong K_*^H(\A(\EE|_U)) \\
                   & \cong K_*^G(\Ind^G_H \A(\EE|_U))  && \cong K_*^G(\A(\EE|_{W})) 
\endaligned$$ and we are finished. \EPf

\appendix
\section{ $\Z_2$-Gradings and Unbounded Multipliers}
In this appendix we collect definitions and results on graded $\cs$-algebras,
tensor products, and unbounded multipliers needed in the text. 
Although some of these results may be found elsewhere ~\cite{BJ83,Bla98,HKT98,Trou00} we include 
them for completeness and to fix notation.

\subsection{Graded $\cs$-algebras}
Recall that a $\cs$-algebra $A$ is {\it $\Z_2$-graded} if it is equipped
with a  $*$-automorphism $\gamma \in \Aut(A)$ such that $\gamma^2 = \id_A$. See
Blackadar~\cite{Bla98} for the basic theory and properties of graded $\cs$-algebras. 

Important examples for us are:
\tightlist
\item The  algebra $\S=C_0(\R)$ of continuous functions on
$\R$, vanishing at infinity, graded according to even and odd functions:
$\gamma(f)(x) = f(-x)$.   
\item The complex Clifford algebra $\Cliff(V)$ of a finite
dimensional
Euclidean vector space $V$ which is the unital $\cs$-algebra generated by a real
subspace $V$ of 
self-adjoint elements  such that  $v^2=\|v\|^2$ for all $v\in V$. The grading
automorphism
is induced by the map $\gamma: V \to V$ given by $\gamma(v) = -v$.
\item Any $\cs$-algebra can be given the {\it trivial} grading
$\gamma = \id_A$. The
$\cs$-algebra $\C$ of complex numbers is always assumed to be trivially graded.
\endtightlist

The {\it graded} commutator of two homogeneous elements $a, b \in A$ is defined 
as $$[a, b] = ab - (-1)^{\deg(a)\deg(b)}ba.$$
We shall say that two elements of $A$  commute if
their graded commutator is zero.

\subsection{Unbounded Multipliers}
Let $A$ be a $\cs$-algebra.  Recall that a {\it multiplier\/} of $A$ is  a 
right $A$-module map $T\colon A\to A$ for which there exists an
`adjoint' right $A$-module map $T^*\colon A\to A$ such that 
$$\langle Tx,y\rangle = \langle x,T^*y\rangle \quad \text{ for all } x,y\in A,$$
where the angle brackets denote the pairing 
$$\langle x,y \rangle=x^*y.$$
Such multipliers form a $*$-algebra under composition, which is denoted
$M(A)$. The operator norm of each multiplier $T$ is finite (by the Closed Graph
Theorem), 
and with this norm $M(A)$ becomes a $\cs$-algebra. It contains $A$ as an
essential ideal.  If $A$ has a $\Z_2$-grading, then $M(A)$ has an induced $\Z_2$-grading by the formula
$$\deg(Tx) \equiv \deg(T) + \deg(x) \ (\bmod 2)$$
for all $T \in M(A)$ and $x \in A$.
See the book~\cite{Ped79} for more details.

We wish to consider some ``unbounded''  multipliers  of $A$.

\begin{dfn}\label{def:usam} $($Compare \cite{BJ83,Con81,HKT98}.$)$ 
An {\it unbounded self-adjoint multiplier\/} of a graded $\cs$-algebra
$A$ is an $A$-linear map $D \colon  \A \to A$ from a dense, right
$\Z_2$-graded $A$-submodule $\A \subset A$ 
into $A$ such that 
\ilist
\item $\langle Dx,y\rangle = \langle x,Dy\rangle  $, for all $x,y\in
 \A$; and 
\item the operators $D\pm iI \colon  \A \to A$  are isomorphisms.
\item $\deg(Dx)\equiv \deg(x)+1 \ (\bmod 2)$, for all $x\in  \A$. 
\endilist
That is, $D$ has grading-degree one with respect to the gradings.
\end{dfn}

\begin{rmk}  The adjective `unbounded'  is conventional, but for 
brevity we shall  often drop this term.  Note that if 
$\A= A$ then the above definition produces a bounded operator.  
\end{rmk}

\subsection{Essentially Self-Adjoint Multipliers}\label{sect:esam}
If $\A$ is a graded, right $A$-submodule of $A$ and if 
$D \colon  \A\to A$  satisfies 
\ilist
\item $\langle Dx,y\rangle = \langle x,Dy\rangle  $, for all $x,y\in
 \A$;
\item $\deg(Dx)\equiv \deg(x)+1 \ (\bmod 2)$, for all $x\in  \A$; 
\item the operators  $D\pm iI$ have dense range;
\endilist 
then the closure $\bar{D}$ of $D$ (in the sense of unbounded operator theory) is
a self-adjoint multiplier.  We shall call $D$ {\it essentially self-adjoint}
and for simplicity we shall usually use the same symbol $D$ for the closure
$\bar{D}$.

Here are some examples.

\noindent $\bullet$ \quad  Denote by $X$ the operator of ``multiplication by
$x$'' on $\R$, so that if $f$ is a function on $\R$ then $Xf(x) = xf(x)$
for all $x \in \R$.  Then $X$ is an unbounded  self-adjoint multiplier of
$\S = C_0(\R)$, whose  domain is the set of all  functions $f\in \S$ such that $Xf\in \S$. 

\noindent $\bullet$ \quad Let $V$ be a Euclidean vector space and let $\CC(V) =
C_0(V, \Cliff_\C V)$
be the $\cs$-algebra of continuous functions from $V$ into the complexified
Clifford algebra
of $V$, with the induced grading. Denote by $C_V : V \to \Cliff_\C(V)$ the
function which
assigns to $v \in V$ the value $v \in \Cliff_\C(V)$. This defines an unbounded
essentially 
self-adjoint multiplier of $\CC(V)$ that is very useful in the formulation of
equivariant 
Bott periodicity for infinite dimensional Euclidean spaces~\cite{HKT98}.

\subsection{Functional Calculus}
If $D$ is an unbounded (essentially) self-adjoint multiplier of $A$ then the
resolvent operators 
$(D\pm iI)^{-1} \colon A \to A$ determine bounded multipliers of
$A$.  If $C_b(\R)$ denotes the $\cs$-algebra of bounded continuous functions on
$\R$ then  there is a unique $*$-homomorphism (called the {\it functional
calculus} of $D$)
$$\gathered C_b(\R) \to M(A)\\ f\mapsto f(D),\endgathered$$ 
mapping the resolvent functions
$(x\pm i)^{-1}$ to $(D\pm iI)^{-1}$.  Compare \cite{Con81}.
This has the property that if $g(x) = xf(x)$ then $g(D) = Df(D)$.  It also
follows from the grading condition (iii) that the above
$*$-homomorphism
is grading preserving, if we grade $C_b(\R)$ by even and odd functions.

\subsection{Graded Tensor Products}
We shall denote by $A\gtimes B$ the {\it maximal\/} graded tensor product of
graded 
$\cs$-algebras. It has the characteristic property that a pair of graded
$*$-homomorphisms $A\to C$ and $B\to C,$
whose images   commute with one another, induces a morphism
$$A\gtimes B\to C.$$ 
In the instances we shall consider, this maximal tensor product is 
also equal to the minimal tensor product, defined with the help of
representations of $A$ and $B$ on graded Hilbert spaces, since one of the
factors will always be nuclear.

If $B$ is {\it evenly} graded, in the sense that there is a self-adjoint unitary
$\varepsilon$ in $ M(B)$ (the multiplier algebra of $B$)  such that for all $b
\in B$,
$$\gamma(b) = \varepsilon b \varepsilon,$$ 
then $A\gtimes B$ is isomorphic to the maximal tensor product of ungraded
$\cs$-algebras, $A\otimes B$,  via the map
$$a\gtimes b \mapsto a\otimes \varepsilon ^{\deg(a)} b.$$ 
This applies, for example, to $B=\K(H)$, where $H$ is a graded Hilbert space.

\subsection{Tensor products and multipliers}
Suppose that $D$ is an (essentially) self-adjoint multiplier of $A$.  
Then the operator $D\gtimes 1$, with domain the algebraic tensor product (over
$\C$)
$\A\gtimes_{alg} B$,  is an
essentially self-adjoint multiplier of $A\gtimes B$.  Note that the functional
calculus
for $D\gtimes 1$ is given by $f(D\gtimes 1)= f(D)\gtimes 1$. 

The following is stated in \cite{HKT98} (see also \cite{BJ83}), but
the proof, which is due to N. Higson, is included below.

\begin{lem}\label{lem:tensormult}
Let $C$ and $D$ be essentially self-adjoint multipliers of $A$ and
$B$. Then $C\gtimes 1 + 1\gtimes D$ is an essentially self-adjoint  
multiplier of $A\gtimes B$, with domain the algebraic tensor product of the
domains
of $C$ and $D$.
\end{lem}

\Pf Let us introduce for each $N > 0$ the bounded multiplier 
$$D_N=D(1+N^{-2}D^2)^{-1/2}$$
(and define $C_N$ similarly).  Note that as $N\to \infty$ the function
$$x_N = x(1+N^{-2}x^2)^{-1/2}$$ converges, uniformly on compact sets, to 
the identity function $x \mapsto x$.  So for $N$ large we expect $D_N$ to
approximate $D$
somehow.  In fact if $\phi$ is a compactly supported function then 
$D_N\phi(D) \to D\phi(D)$ in norm.  

Let us make a further observation: for every $y\in A\gtimes B$ there 
is some compactly supported $\phi$ so that 
$$y\approx \phi(C)\gtimes \phi(D) y$$
(in other words  given $\varepsilon >0$ we can find $\phi$ so that the two
sides are
within $\varepsilon$ in norm).  

Now given $y$ choose $\phi$ as above then choose $N$ so that 
$D_N\phi(D)\approx  D\phi(D)$.  By simple
$\cs$-algebra theory  the bounded multiplier  $C_N\gtimes 1 + 1\gtimes
D_N$ is
essentially  self-adjoint in the sense of 1.4.  So  there 
is some $x\in \A\gtimes_{alg} \B$ with 
$$ y\approx (iI+C_N\gtimes 1 + 1\gtimes D_N)x.$$
But then 
 $$ \aligned y \approx \phi(C)\gtimes \phi(D) y&\approx \phi(C)\gtimes
\phi(D)(iI+C_N\gtimes 1 + 1\gtimes D_N)x\\
&=(iI+C_N\gtimes 1 + 1\gtimes D_N) \phi(C)\gtimes
\phi(D)x\\
&\approx (iI+C\gtimes 1 + 1\gtimes D) \phi(C)\gtimes
\phi(D)x,\endaligned $$
so that $y$ is in the closure of the range of $(iI+C\gtimes 1 + 1\gtimes
D)$ as
required. \EPf

\begin{cor}\label{cor:exp}
If $C$ and $D$ are essentially self-adjoint multipliers of $A$ and
$B$, then
$$\exp(-(C\gtimes 1 + 1\gtimes D)^2) = \exp(-C^2)\gtimes \exp(-
D^2).$$
\end{cor}

\subsection{Graded $C_0(M)$-algebras}

Let $M$ be a locally compact space. We assume that $C_0(M)$ has the trivial grading.

\begin{dfn}\label{def:grddCo(M)alg} A $\cs$-algebra $A$ is called a {\bf $\Z_2$-graded $C_0(M)$-algebra}
if $A$ is $\Z_2$-graded and there is a $*$-homomorphism $C_0(M) \to ZM(A)$ 
of degree zero which is nondegenerate in the sense that $C_0(M) A$ is dense in $A$.
\end{dfn}

We now introduce a new tensor product, which is a graded version of the balanced tensor product considered
by others \cite{Blan95,EW97,RW98}.

\begin{dfn}\label{def:gtimes_M}
Let $A$ and $B$ be graded $C_0(M)$-algebras. We define the balanced $($graded$)$ tensor product $A \gtimes_{_M} B$
to be the quotient of $A \gtimes B$ by the $\Z_2$-graded ideal $J_M$ generated by
$$\{(f \cdot a) \gtimes b - a \gtimes (f\cdot b) | f \in C_0(M), a \in A, \text{ and } b \in B\}.$$
There is then a short exact sequence
$$\xymatrix{ 0 \ar[r] & J_M \ar[r] & A \gtimes B \ar[r]^-{q_{_M}} & A \gtimes_{_M} B \ar[r] & 0 }$$
of $\Z_2$-graded $\cs$-algebras.
\end{dfn}

\noindent There is an isomorphism $C_0(M) \gtimes_{M} A \cong A$ given by the map $f \gtimes a \mapsto f \cdot a$.

This tensor product has the universal property that if $\phi : A \to C$ and $\psi : B \to C$
are graded $C_0(M)$-linear $*$-homomorphisms with the property that their ranges {\it graded} commute, 
then there is a unique graded $*$-homomorphism
$$\phi \gtimes_{_M} \psi : A \gtimes_{_M} B \to M(C) : a \gtimes_{M} b \mapsto \phi(a)\psi(b)$$
where $a \gtimes_{_M} b$ denotes the image of $a \gtimes b$ in $A \gtimes_{_M} B$.

\begin{prop}\label{prop:balanced}
Let $M$ be a Riemannian manifold. If $E$ and $F$ are finite rank Euclidean bundles on $M$,
then there is a $\Z_2$-graded $C_0(M)$-algebra isomorphism
$$C_0(M, \Cliff(E \oplus F)) \cong C_0(M, \Cliff(E)) \gtimes_{_M} C_0(M, \Cliff(F)).$$
\end{prop}

\Pf There is an $\Z_2$-graded isomorphism of Clifford algebra bundles
$$\Cliff(E \oplus F) \cong \Cliff(E) \gtimes \Cliff(F)$$
which follows by Proposition 1.5 \cite{LM89}. The result now
follows by applying the same argument as used in Proposition 2.6 \cite{GVF01}
and the universal property of the balanced tensor product. \EPf

\begin{dfn} Let $p : E \to M$ be a continuous map of locally compact spaces.
If $A$ is a $\Z_2$-graded $C_0(M)$-algebra, then we define the {\bf pullback
$\cs$-algebra} $p^*A$ by the formula
$$p^*A = C_0(E) \gtimes_M A$$
where $C_0(E)$ is given a $C_0(M)$-algebra structure via $p^* : C_0(M) \to C_b(E) \cong M(C_0(E))$.
Note that since $C_0(M)$ is trivially graded, this is the same definition as
given in \cite{EW97, RW98}.
\end{dfn}

The relationship to pullback vector bundles is contained in the next result.

\begin{prop}\label{prop:pullbackalg}
Let $E$ be a smooth finite rank Euclidean bundle on the Riemannian manifold $M$.
If $p : N \to M$ is a smooth map of Riemannian manifolds, there is an 
isomorphism
$$p^*C_0(M, \Cliff(E)) \cong C_0(N, \Cliff(p^*E))$$
of $\Z_2$-graded $\cs$-algebras.
\end{prop}

\Pf By universality, there is a pullback isomorphism of
Clifford bundles \break $\Cliff(p^*E) \cong p^*\Cliff(E)$.
The result now follows from Proposition 2.12~\cite{GVF01}. \EPf

\subsection{Pushing Forward Multipliers} 
 Suppose we are given a $*$-homomorphism
$\phi : A'\to M(A)$ with the property that
$A'A$ is dense in
$A$. Then there is an induced $*$-~homomorphism 
$M(A')\to M(A)$ extending this, characterized by $T(a'a)= T(a')a$. 
See~\cite{Ped79}.
Similarly, an essentially self-adjoint multiplier $D'$ of $A'$, with domain 
$\A'\subset A'$, determines an essentially self-adjoint operator on  $A$, 
which we shall denote $D$ here, by the formula
$$\qquad \qquad D(a'a)= (D'a')a  \qquad (\text{where $a'\in  \A' $ and $\ 
a\in A$}).$$  
Its domain is $\A' A$. The functional calculus maps for the two operators
$D'$ and $D$ 
(or, to be more precise, for the  closures of $D'$ and $D$) are related by
the following commutative diagram:
$$\CD C_b(\R) @>f\mapsto f(D')>> M(A') \\
@V=VV @VVV\\
C_b(\R) @>>f\mapsto f(D )> M(A )
\endCD.$$


\bibliographystyle{amsplain}
\providecommand{\bysame}{\leavevmode\hbox to3em{\hrulefill}\thinspace}

\end{document}